\documentclass[11pt]{article}
\usepackage{latexsym,amsmath,amssymb,amsfonts,amsbsy,fancyhdr}
\usepackage{theorem}
\usepackage{amscd}
\usepackage{enumerate}
\usepackage{mathrsfs}
\usepackage{dsfont}
\usepackage[active]{srcltx}
\usepackage[colorlinks=true, pdfstartview=FitV, linkcolor=black, citecolor=black, urlcolor=black]{hyperref}
\usepackage[dvips]{graphicx}
\usepackage{epsf}
\usepackage{nicefrac}
\usepackage{marvosym}
\usepackage{subfigure}

\newcommand{\R}{{\mathbb R}}
\newcommand{\rp}{\R^+}
\newcommand{\rpn}{\R^+_0}

\newcommand{\N}{{\mathbb N}}

\renewcommand{\d}{\text{ d}}
\newcommand{\dd}{\text{d}}

\newtheorem{theorem}{Theorem}
\newtheorem{proposition}{Proposition}
\newtheorem{corollary}[proposition]{Corollary}
\newtheorem{lemma}[proposition]{Lemma}

\theorembodyfont{\rm}
\newtheorem{definition}[proposition]{Definition}
\newtheorem{example}[proposition]{Example}

\newtheorem{remark}[proposition]{Remark}
\newenvironment{rem}{\begin{remark}}{\hfill$\lozenge$\end{remark}}
\newtheorem{proof}{Proof}

\newenvironment{pf}{\begin{proof}}{\hfill$\square$\end{proof}}

\allowdisplaybreaks

\DeclareMathAlphabet{\mathcal}{OMS}{cmsy}{m}{n}

\begin{document}

\title{\Large \bf One-sided FKPP travelling waves in the context of homogeneous fragmentation processes}
\author{Robert Knobloch\thanks{Department of Mathematics, Saarland University, PO Box 151150, 66041 Saarbr\"ucken, Germany
\newline
e-mail: knobloch@math.uni-sb.de} }
\date{\today}
\maketitle

\begin{abstract}
In this paper we introduce the one-sided FKPP equation in the context of homogeneous fragmentation processes. The main result of the present paper is concerned with the existence and uniqueness of one-sided FKPP travelling waves in this setting. Moreover, we prove some analytic properties of such travelling waves. Our techniques make use of  fragmentation processes with killing, an associated product martingale as well as various properties of 
L\'evy processes.
\end{abstract}




\section{Introduction}
This paper deals with an integro-differential  equation that is defined in terms of the dislocation measure of a fragmentation process. Given its probabilistic interpretation we consider this equation as an analogue of the one-sided Fisher-Kolmogorov-Petrovskii-Piscounov (FKPP) travelling wave equation in the context of  fragmentation processes.
In particular, we are concerned with the existence and uniqueness of solutions to this equation in the setting of conservative or dissipative homogeneous fragmentation processes and we derive certain analytic properties of the solutions when they exist.

The FKPP travelling wave equation in the context of  fragmentations
has a similar probabilistic interpretation as the classical FKPP travelling wave equation whose probabilistic interpretation is related to branching Brownian motions, see Section~\ref{s.cFKPPe}. In this respect we also refer to  \cite{104}, where the two-sided FKPP travelling wave equation for conservative homogeneous fragmentations is studied.

It turns out that  solutions in the setting of fragmentation processes have similar properties as their classical counterparts. However, the techniques of proving existence and uniqueness results are very different between the two cases, not least because the equations 
differ significantly. Indeed, whereas the classical FKPP travelling wave equation is a differential equation of second order, the FKPP travelling wave equation in our setting is an integro-differential equation of first order. This difference results from the non-diffusive behaviour of fragmentation processes and the more complicated jump structure of fragmentations in comparison with  branching Brownian motions. 

In the context of  homogeneous fragmentation processes we prove the existence and uniqueness of one-sided travelling waves  within a certain range of wave speeds. More precisely, the problem we are concerned with in this paper can be roughly described as follows. 
Consider the integro-differential equation
\[
cf'(x)+\int_{\mathcal P}\left(\prod_{n}f(x+\ln(|\pi_n|))-f(x)\right)\mu(\dd\pi)=0
\]
for certain $c\in\rp:=(0,\infty)$ and all $x\in\rpn:=[0,\infty)$, where the product is taken over all $n\in\N$ with $|\pi_n|\in\rp$. Here the space $\mathcal P$ is the space of partitions $(\pi_n)_{n\in\N}$ of $\N$ and $\mu$ is the so-called dislocation measure on $\mathcal P$. This notation is introduced in more detail in the next section. We are interested in solutions $f:\R\to[0,1]$ of the above equation that satisfy
\[
f|_{\rpn}\in C^1(\rpn,[0,1])\qquad\text{and}\qquad f|_{(-\infty,0)}\equiv1
\]
as well as the boundary condition
\[
\lim_{x\to\infty}f(x)=0.
\] 
Roughly speaking, the main result of this paper states that there is some constant $c_0>0$ such that there exists a unique solution of the above boundary value problem for every $c>c_0$ and there does not exist such a solution for any $c\le c_0$. Our approach is based on using fragmentation processes with killing at an exponential barrier. These processes have been studied in \cite{KK12} and we briefly describe the corresponding concepts below.

The outline of this paper is as follows. In the next section we give a brief introduction to homogeneous fragmentation processes as well as appropriately killed versions of such fragmentations.  Subsequently, in Section~\ref{s.FKPPfrag} we introduce the one-sided FKPP equation in our setting and we state our main results. Afterwards, in the fourth section we provide some motivation for the problems considered in this paper by explaining some related results that are known in the literature on fragmentation processes and branching Brownian motions, respectively. In Section~\ref{s.fa} we show how the existence and uniqueness of one-sided FKPP travelling waves for fragmentation processes can be obtained if the dislocation measure is finite. This provides some motivation for the existence and uniqueness result in the setting of general fragmentation processes, which this paper is mainly concerned with. The subsequent three sections are devoted to the proofs of our main results.


Throughout the present paper we adopt the notation $\R_\infty:=[-\infty,\infty)$ as well as the conventions $\ln(0):=-\infty$ and $\inf(\emptyset):=\infty$. The notation $C^n$, $n\in\N_0:=\N\cup\{0\}$, refers to the set of $n$-times continuously differentiable functions. The integral of a real-valued function $f$ with respect to the Lebesgue measure on a set $[s,t]\subseteq\R$ is denoted by $\int_{[s,t]}f(u)\dd u$ and $\int_s^tf(u)\dd u$ denotes the Riemann integral. The operators $\land$ and $\lor$ refer to the minimum and maximum, respectively. Furthermore, we shall use the abbreviation DCT for the  dominated convergence theorem. All the  random objects are assumed to be defined on a complete probability space $(\Omega,\mathscr F,\mathbb P)$.

\section{Homogeneous fragmentation processes with killing}\label{s.khfp}
In this section we provide a brief introduction to partition-valued fragmentation processes and we present the main tools that we need in the subsequent sections. In addition, we introduce a specific killing mechanism for these processes. The advantage of partition-valued fragmentation processes compared to so-called mass fragmentations is their explicit  genealogical structure of blocks. This structure is crucial for the killing mechanism that we  introduce below.

Regarding the state space of partition-valued fragmentation processes let $\mathcal{P}$ be the space of partitions $\pi=(\pi_n)_{n\in\N}$ of $\N$, where the blocks of $\pi$ are ordered by their least element such that $\inf(\pi_i)<\inf(\pi_j)$ if $i<j$. For every $\pi\in\mathcal P$ let $(|\pi|^\downarrow_n)_{n\in\N}$ be the decreasing reordering of the sequence given by
\[
|\pi_n|:=\limsup_{k\to\infty}\frac{\sharp(\pi_n\cap\{1,\ldots,k\})}{k}
\]
for every $n\in\N$, where $\sharp$ denotes the counting measure on $\N$.

Throughout this paper we consider a homogeneous $\mathcal P$-valued fragmentation process $\Pi:= (\Pi(t))_{t\in\rpn}$, where $\Pi(t) = (\Pi_n(t))_{n\in\N}$, and we denote by $\mathscr F:=(\mathscr F_t)_{t\in\rpn}$ the completion of the filtration generated by $\Pi$.  Homogeneous  $\mathcal{P}$-valued fragmentations are exchangeable Markov processes that were introduced in  \cite{85}, see also \cite{92}.   Bertoin showed in \cite{85} that the distribution of $\Pi$ is determined by some constant $d\in\rpn$ (the {\it rate of erosion} which describes the drift of $\Pi$) and a $\sigma$-finite measure $\nu$ (the so-called {\it dislocation measure} that indirectly describes the jumps of $\Pi$) on the infinite simplex
 \[
\mathcal S:=\left\{{\bf s}:=(s_n)_{n\in\N}:s_1\ge s_2\ldots\ge0,\,\sum_{n\in\N}s_n\le1\right\},
 \] 
such that $\nu(\{(1,0,\ldots)\})=0$ and
\begin{equation}\label{e.levymeasure}
\int_{\mathcal{S}}(1-s_1)\nu(d{\bf s})<\infty.
\end{equation}

The process $\Pi$ is said to be {\it conservative} if $\nu(\sum_{n\in\N}s_n<1)=0$, i.e. if there is no loss of mass by sudden dislocations, and {\it dissipative} otherwise. In this paper we allow for both of these cases. 

{\it Throughout this paper we assume that $d=0$  as well as $\nu({\bf s}\in\mathcal S:s_2=0)=0$.} 
\\
In view of the forthcoming assumption (\ref{e.L.1.Bertoin}) this enables us to resort to the results of \cite{KK12}, where the same assumptions are made. Let us mention that the assumption $d=0$ does not result in any loss of generality, see Remark~\ref{r.drift}.

Consider the  exchangeable partition measure $\mu$ on $\mathcal{P}$ given by
\[
\mu(d\pi) = \int_{\mathcal{S}}\varrho_{\bf s}(d\pi)\nu(d{\bf s}),
\]
where $\varrho_{\bf s}$ is the law of Kingman's paint-box based on ${\bf s}\in\mathcal S$. Similarly to $\nu$ the measure $\mu$ describes the jumps of $\Pi$, although more directly,  and is also referred to as {\it dislocation measure}. In \cite{85} Bertoin showed that the homogeneous  fragmentation process $\Pi$ is characterised by a Poisson point process. More precisely, there exists a $\mathcal P\times\mathbb N$-valued Poisson point process $(\pi(t),\kappa(t))_{t\in\rpn}$\label{p.pi_t} with characteristic measure $\mu\otimes\sharp$  such that $\Pi$ changes state only at the times $t\in\rpn$ for which an atom $(\pi(t),\kappa(t))$ occurs in $(\mathcal P\setminus(\N,\emptyset,\ldots))\times\mathbb N$. At such a time $t\in\rpn$ the sequence $\Pi(t)$ is obtained from $\Pi(t-)$ by replacing its $\kappa(t)$-th term, $\Pi_{\kappa(t)}(t-)\subseteq\N$, with the restricted partition $\pi(t)|_{\Pi_{\kappa(t)}(t-)}$ and reordering the terms such that the resulting
partition of $\N$ is an element of $\mathcal P$. We denote the possible random jump times of $\Pi$, i.e. the times at which the abovementioned Poisson point process has an atom in $(\mathcal P\setminus(\N,\emptyset,\ldots))\times\mathbb N$, by $(t_i)_{i\in\mathcal I}$, where the index set $\mathcal I\subseteq\rpn$ is countable. 

Moreover, by exchangeability, the limit
\[
|\Pi_n(t)|:=\lim_{k\to\infty}\frac{\sharp(\Pi_n(t)\cap\{1,\ldots,k\})}{k},
\]
referred to as {\it asymptotic frequency}, exists $\mathbb P$-a.s. simultaneously for all $t\in\rpn$ and $n\in\N$. Let us point out that the concept of asymptotic frequencies provides us with a notion of {\it size} for the blocks of a $\mathcal P$-valued fragmentation process. In Theorem~3 of \cite{85} Bertoin showed that the process $(-\ln(|\Pi_1(t)|))_{t\in\rpn}$ is a killed subordinator, a fact we shall make use of below.   
 
\begin{figure}[htb]
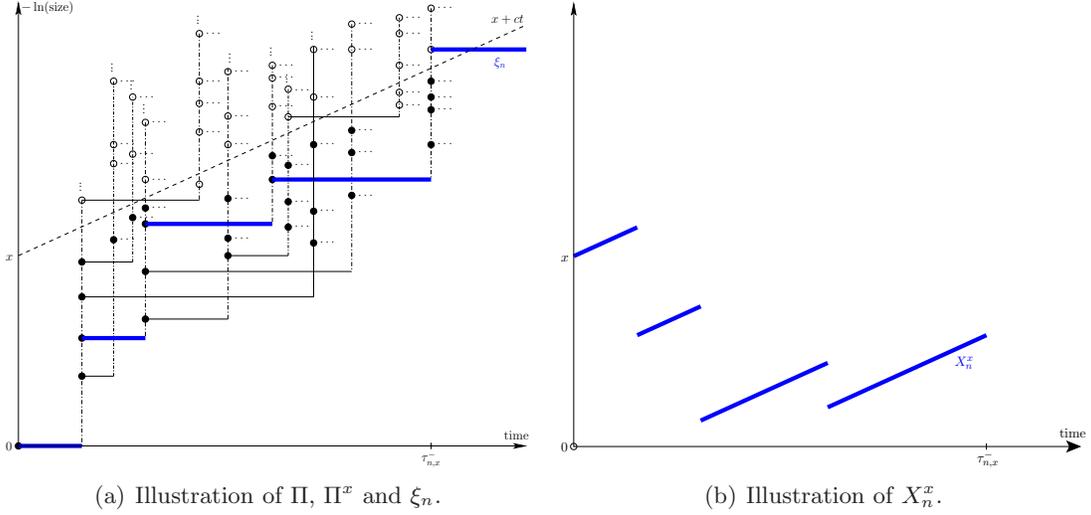
\centering
\subfigure[Illustration of $\Pi$, $\Pi^x$ and $\xi_n$.]{\resizebox{7cm}{!}{\input{pic4long_125_c.pstex_t}}\label{f.kfp.2a}}\quad
\subfigure[Illustration of $X^x_n$.]{\resizebox{7cm}{!}{\input{pic4bkilled_125c.pstex_t}}\label{f.kfp.3killed}}\quad
\caption[Spectrally negative L\'evy process associated with a fragmentation]{Illustration (a) depicts $\Pi$, $\Pi^x$ and $\xi_n$ for a finite dislocation measure. The black dots indicate the blocks that belong to $\Pi^x$, whereas $\Pi$ additionally contains the white dots. The dashed line constitutes the killing line with slope $c>0$ starting at $x\in\rpn$. In addition, the bold piecewise constant graph illustrates the subordinator $\xi_n$. Notice that $\tau^-_{n,x}$ is the first time that $\xi_n$ is above the killing line, i.e. it is the killing time of the block containing $n\in\N$. Illustration (b) shows the shifted and killed spectrally negative L\'evy process $X^x_n$. After time $\tau^-_{n,x}$ the process $X^x_n$ takes the value $-\infty$ and thus does not appear in the illustration anymore.}
\end{figure}

For the time being, let $x\in\R_\infty$. In this paper we are concerned with a specific procedure of killing blocks of $\Pi$, see Figure~\ref{f.kfp.2a}, that was introduced in \cite{KK12}. More precisely, for $c>0$
a block $\Pi_n(t)$ is killed, with cemetery state $\emptyset$,  at the moment of its creation $t\in\rpn$ if $|\Pi_n(t)|<e^{-(x+ct)}$. We denote the resulting fragmentation process with killing by $\Pi^x:=(\Pi^x(t))_{t\in\rpn}$ and the cemetery state  of $\Pi^x$ is  $(\emptyset,\ldots)$.
Note that possibly $\Pi^x(t)\not\in\mathcal P$ as $\bigcup_{n\in\N}\Pi^x_n(t)\subsetneq\N$ is possible due to the killing of blocks. 

We denote by $\zeta^x$, $x\in\R_\infty$, the random extinction time of $\Pi^x$, i.e. $\zeta^x$ is the supremum of all the killing times of individual blocks. The question whether $\zeta^x$ is finite or infinite was considered in Theorem~2 of \cite{KK12}, see also Proposion~\ref{positivesurviaval} below. Furthermore, define a function $\varphi:\R\to[0,1]$ by
\begin{equation}\label{e.phi}
\varphi(x):=\mathbb P(\zeta^x<\infty)
\end{equation}
for all $x\in\R_\infty$.
The function $\varphi$ will be of utmost interest in the present paper. Let us point out that $\varphi$ depends on the drift $c>0$ of the killing line. However, in order to keep the notation as simple as possible, we omit this dependence in the notation as the constant $c$ does not vary within results or proofs. Note that if $x<0$, then $\zeta^x=0$ and thus $\varphi(x)=1$. 

Let us remark that we could choose a non-zero rate $d$ of erosion by changing the slope $c$ of the killing line:

\begin{rem}\label{r.drift}
The  results of this paper remain valid if we omit the assumption $d=0$ and replace the slope $c>0$ by $c_d:=c+d$. Consequently, the assumption $d=0$ is merely made  for the sake of simplicity, but does not restrict the generality of our results.
\end{rem}

Set
\[
\underline p:=\inf\left\{p\in\R:\int_{\mathcal S}\left|1-\sum_{n\in\N}s^{1+p}_n\right|\nu(\dd{\bf s})<\infty\right\}\in[-1,0]
\]
and for any $p>\underline p$ define
\[
\Phi(p):=\int_{\mathcal S}\left(1-\sum_{n\in\N}s^{1+p}_n\right)\nu(\dd{\bf s})
\]
as well as
\[
\Phi(\underline p):=\lim_{p\downarrow \underline p}\Phi(p).
\]
Moreover, for each $p\in[\underline p,\infty)$ set
\begin{equation}\label{c_p.2}
c_p:=\frac{\Phi(p)}{1+p}.
\end{equation}
Throughout this paper we assume that there exists some $p\in(\underline p,\infty)$ such that
\begin{equation}\label{e.L.1.Bertoin}
(1+p)\Phi'(p)>\Phi(p),
\end{equation}
where $\Phi'$ denotes the derivative of $\Phi$. Let us point out that a sufficient condition for (\ref{e.L.1.Bertoin}) is the existence of some $p^*\in[\underline p,\infty)$ such that $\Phi(p^*)=0$. In particular,  (\ref{e.L.1.Bertoin}) holds if $\Pi$ is conservative.
In view of (\ref{e.L.1.Bertoin}) the same line of argument as in Lemma~1 of \cite{89} yields the existence of  a unique solution of the equation 
\begin{equation}\label{e.bar_p}
(1+p)\Phi'(p)=\Phi(p)
\end{equation} 
on $(\underline p,\infty)$. We denote this unique solution of (\ref{e.bar_p}) by $\bar p$. The definition in (\ref{c_p.2}) then entails that $c_{\bar p}=\Phi'(\bar p)$. According to \cite{KK12} the  fragmentation process with killing survives with positive probability if the drift of the killing line is greater than $c_{\bar p}$ and becomes extinct almost surely otherwise. 
\begin{proposition}[Theorem~2 of \cite{KK12}]\label{positivesurviaval}
If $c>c_{\bar p}$, then $\varphi(x)\in(0,1)$ for all $x\in\rpn$. If, on the other hand, $c\le c_{\bar p}$, then $\varphi\equiv1$.
\end{proposition}

For any  $t\in\rpn$ we denote by  $B_n(t)$ the block of $\Pi(t)$ that contains the element $n\in\N$. According to Theorem~3 (ii) in \cite{85} it follows by means of the exchangeability of $\Pi$ that under $\mathbb P$ the process 
\[
\xi_n:=(-\ln(|B_n(t)|)
)_{t\in\rpn}\,,
\] 
cf. Figure~\ref{f.kfp.2a}, is a killed subordinator with Laplace exponent $\Phi$, cemetery state $\infty$ and killing rate
\[
\int_{\mathcal S}\left(1-\sum_{k\in\N}s_k\right)\nu(\dd{\bf s}).
\] 
Hence, the process $X_n:=(X_n(t))_{t\in\rpn}$, defined by 
\[
X_n(t):=ct-\xi_n(t)
\] 
for all $t\in\rpn$, is a spectrally negative L\'evy process  of bounded variation. Let $\mathcal I_n\subset\mathcal I$ be such that the jump times of $X_n$ are given by $(t_i)_{i\in\mathcal I_n}$. Note that $(t_i)_{i\in\mathcal I_n}$ are precisely the times when the subordinator $\xi_n$ jumps. For $n\in\N$ and $x\in\rpn$ we shall make use of the shifted and killed process $X^x_n:=(X^x_n(t))_{t\in\rpn}$, see Figure~\ref{f.kfp.3killed}, given by
\[
X^x_n(t):=(X_n(t)+x)\mathds1_{\{\tau^-_{n,x}>t\}}-\infty\cdot\mathds1_{\{\tau^-_{n,x}\le t\}}=\left(x+ct+\ln(|B_n(t)|)\right)\mathds1_{\{\tau^-_{n,x}>t\}}-\infty\cdot\mathds1_{\{\tau^-_{n,x}\le t\}}
\]
for each $t\in\rpn$, where 
\[
\tau^-_{n,x}:=\inf\{t\in\rpn:X_n(t)<-x\}
\qquad\text{as well as}\qquad
\infty\cdot0:=0.
\]

For every $t\in\rpn$ set 
\[
\mathcal N^x_t:=\left\{n\in\N:\left[t<\tau^-_{n,x}\right]\land\left[\exists\,k\in\N:n=\min\Pi^x_k(t)\right]\right\}.\label{p.Nxt}
\]
That is to say,  $\mathcal N^x_t$ consists of all the indices of blocks $B_n(t)$ that are not yet killed by time $t$. Let us remark that the first condition ``$t<\tau^-_{n,x}$''  ensures that  the block containing $n\in\N$ is still alive at time $t$ and the second condition ``$\exists\,k\in\N:n=\min\left(\Pi^x_k(t)\right)$'' is used to avoid considering the same block multiple times. More precisely, for a block $B_n(t)$ that is alive at time $t\in\rpn$ only its least element is an element of $\mathcal N^x_t$. Without this condition all elements of $B_n(t)$ would be in $\mathcal N^x_t$. 
Note that in \cite{KK12} the notation $\mathcal N^x_t$ refers to a different ordering of the same set of blocks and the set of indices that we denote here by $\mathcal N^x_t$ is the quotient space $\widetilde{\mathcal{N}}_t^x/\sim$ in \cite{KK12}.

Throughout this paper we shall repeatedly need an estimate regarding the number of fragments alive at 
a given time $t\in\rpn$. To this end, set $N^x_t:=\text{card}(\mathcal N^x_t)$ and 
observe that $N^x_t<\infty$ $\mathbb P$-almost surely. Indeed, since $\sum_{n\in\N}|\Pi_n(t)|\le1$ we infer that $|\Pi_n(t)|\ge e^{-(x+ct)}$ for at most $e^{x+ct}$-many $n\in\N$. Hence,
\begin{equation}\label{e.N_x}
N^x_t\le e^{x+ct}. 
\end{equation}

\section{The one-sided FKPP equation for  fragmentations}\label{s.FKPPfrag}
In this section we establish the set-up for our considerations by defining the FKPP travelling wave  equation in the context of fragmentation processes. Moreover, this section is devoted to presenting our main results. The main problem addressed in this paper is to find a range of wave speeds for which we can prove the existence of a unique travelling wave solution to the one-sided FKPP equation for homogeneous fragmentations as defined below. In order to tackle this problem we shall derive a connection between solutions of the abovementioned FKPP equation and a product martingale that was introduced in \cite{KK12}. Furthermore, we aim at studying some analytic properties of such travelling waves.
\subsection{Set-up}
For any function $f$ on some subset of $\R_\infty$ set
\[
\mathcal C_f:=\left\{x\in\rp:f'(x)\text{ exists}\right\},
\]
where $f'$ denotes the derivative of $f$.
Since we do not know a priori whether the functions $f$ we are interested in are differentiable, we need to define the integro-differential equations below for arguments in $\mathcal C_f$. Regarding solutions of such an equation  we shall be particularly concerned with the function $\varphi$, given by  (\ref{e.phi}), and our main results show in particular  that  $\mathcal C_\varphi=\rp$.

For functions $u:\rpn\times\R_\infty\to[0,1]$, with $u(t,\cdot)|_{[-\infty,ct)}\equiv1$ for each $t\in\rpn$ and $u(0,\cdot)|_{\rpn}=g$ for some continuous  function $g:\rpn\to[0,1]$, consider the  initial value problem
\begin{equation}\label{e.FKPP.0}
\frac{\partial u}{\partial t}(t,x)=\int_{\mathcal P}\left(\prod_{n\in\N}u(t,x+\ln(|\pi_n|))-u(t,x)\right)\mu(\dd\pi) 
\end{equation}
for all  $x\in\rpn$ and $t\in\mathcal C_{u(\cdot,x)}$. We call this initial value problem \emph{one-sided FKPP equation for  fragmentation processes}. Here we are interested in the so-called \emph{FKPP travelling wave solutions} of (\ref{e.FKPP.0}) with wave speed $c\in\rpn$, that is in solutions of (\ref{e.FKPP.0}) which are of the form $u(t,x)=f(x-ct)$ for all $t,x\in\rpn$ with $x-ct\ge0$.
\begin{definition}
A \emph{one-sided  FKPP travelling wave for fragmentation processes} is a monotone function $f:\R_\infty\to[0,1]$, with $f|_{[-\infty,0)}\equiv1$, that satisfies the following \emph{one-sided FKPP travelling wave equation} 
\begin{equation}\label{e.FKPP.1}
cf'(x)+\int_{\mathcal P}\left(\prod_{n\in\N}f(x+\ln(|\pi_n|))-f(x)\right)\mu(\dd\pi)=0
\end{equation}
for all $x\in\mathcal C_f$ with the boundary condition
\begin{equation}\label{e.FKPP.2}
\lim_{x\to\infty}f(x)=0.
\end{equation}
\end{definition}
Let us now introduce an operator whose definition  is inspired by the integro-differential equation~(\ref{e.FKPP.1}). To this end, let $\mathcal D_L$ be the set of all functions $f:\R_\infty\to[0,1]$, with $f|_{[-\infty,0)}\equiv1$, for which the mapping
\[
\pi\mapsto\prod_{n\in\N}f(x+\ln(|\pi_n|)-f(x)
\]
on $\mathcal P$ is $\mu$-integrable. Then we define an integral operator $L$ with domain $\mathcal D_L$ by
\[
Lf(x):=\int_{\mathcal P}\left(\prod_{n\in\N}f(x+\ln(|\pi_n|)-f(x)\right)\mu(\dd\pi)\label{p.L}
\]
for each $f\in\mathcal D_L$ and all $x\in\rpn$. 

Recall that the upper Dini derivative (from above) of a function $f:\R_\infty\to[0,1]$ is defined by
\[
f'_+(x):=\limsup_{h\downarrow 0}\frac{f(x+h)-f(x)}{h}
\]
for all $x\in\R$. On this note, observe that the Dini derivative is well defined, but may take the value $\infty$ or $-\infty$.

The following class of monotone functions plays a crucial role in the analysis of the one-sided FKPP travelling wave equation.
\begin{definition}\label{d.T}
We denote by $\mathcal T$\label{p.Dp} the set consisting of all nonincreasing functions $f:\R_\infty\to[0,1]$, with $f|_{[-\infty,0)}\equiv1$ and such that $f|_{\rp}$ is continuous, that satisfy (\ref{e.FKPP.2}) as well as
\begin{equation}\label{finite_derivative}
\sup_{x\in[s,t]}|f'_+(x)|<\infty
\end{equation}
for all $s,t\in\rp$.
\end{definition}

\subsection{Main results}
Here we present the main results of this paper. For this purpose, let us consider the following process which several of our proofs make use of: For any function $f:\R_\infty\to[0,1]$ and $x\in\rpn$ let $Z^{x,f}:=(Z^{x,f}_t)_{t\in\rpn}$ be given by
\begin{equation}\label{e.prod_martingale}
Z^{x,f}_t:=\prod_{n\in\mathcal N^x_t}f\left(X^x_n(t)\right)
\end{equation}
for all $t\in\rpn$. This process was considered in Section~5 of \cite{KK12} and in some proofs we shall resort to Theorem~10 in \cite{KK12}, which is concerned with the martingale property of $Z^{x,f}$. In this spirit, the first main result of the present paper reads as follows:
\begin{theorem}\label{p.t.1.2}
Let $c>c_{\bar p}$. In addition, let $f\in\mathcal T$ and assume that $Z^{x,f}$ is a martingale. Then  $f$ solves (\ref{e.FKPP.1}).
\end{theorem}

The above theorem will be proven in Section~\ref{s.cspm}. The following result, whose proof will be provided in Section~\ref{s.ap}, deals with some analytic properties of one-sided FKPP travelling waves.

\begin{theorem}\label{t.ap}
Every one-sided FKPP travelling wave $f\in\mathcal T$ is right-continuous at $0$ and the function $f|_{\rp}$ is strictly monotonically decreasing and continuously differentiable.
\end{theorem} 

In particular, it follows from Theorem~\ref{t.ap} that $-cf'=Lf$ on $\rp$ for every one-sided FKPP travelling wave $f\in\mathcal T$, where $c>0$ denotes the wave speed. 

The main goal of this paper is to establish the existence of a unique  travelling wave in $\mathcal T$ to (\ref{e.FKPP.0}) with wave speed $c$ for $c>c_{\bar p}$ as well as the nonexistence of such a travelling wave with wave speed $c\le c_{\bar p}$. More specifically, the following result states that the extinction probability of the  fragmentation process with killing solves equation~(\ref{e.FKPP.1}) with boundary condition~(\ref{e.FKPP.2}) for $c>c_{\bar p}$ and, moreover, is the only such function. Recall the function $\varphi$ defined in  (\ref{e.phi}).
\begin{theorem}\label{t.mainresult.1}
If $c>c_{\bar p}$, then there exists a unique one-sided FKPP travelling wave in $\mathcal T$ with wave speed $c$, given by  $\varphi$. 
On the other hand, if $c\le c_{\bar p}$, then there is no one-sided FKPP travelling wave in $\mathcal T$ with wave speed $c$.
\end{theorem}
We shall prove Theorem~\ref{t.mainresult.1} in Section~\ref{s.mainresult.1}.  In view of the forthcoming Remark~\ref{r.BHK09} this theorem shows that one-sided FKPP travelling waves exist precisely for those positive wave speeds 
for which there do not exist two-sided travelling waves. Moreover, in view of \cite[Theorem~4]{KK12}, Theorem~\ref{t.mainresult.1} shows that travelling wave solutions exist exactly for those wave speeds that are larger than the asymptotic decay of the largest fragment in the  fragmentation process with killing on the event of survival. Let us further point out that if $\Pi$ is conservative, i.e. if loss of mass by sudden dislocations is excluded, then in view of (3) in \cite{HKK10} equation (\ref{e.FKPP.1}) can be written as
\[
cf'(x)+\int_{\mathcal S}\left(\prod_{n\in\N}f(x+\ln(s_n))-f(x)\right)\nu(\dd{\bf s})=0,
\]
which is the analogue of the two-sided FKPP equation considered in \cite{104} in the conservative setting, cf. (\ref{e.FKPP.0.a}).

\section{Motivation -- The classical FKPP equation}\label{s.cFKPPe}
In order to present the framework in which Theorem~\ref{t.mainresult.1} should be seen let us now briefly mention some known results that are related to our work. To this end we denote by $C^{1,2}(\rpn\times A,[0,1])$\label{p.C^1_2}, $A\subseteq\R$, the space of all functions $f:\rpn\times A\to[0,1]$ such that $f(x,\cdot)\in C^2(A,[0,1])$ and $f(\cdot,y)\in C^1(\rpn,[0,1])$ for all $x\in\rpn$ and $y\in A$. 

The classical FKPP equation 
in the form that is of most interest for us, cf. \cite{McK75}, is the following 
parabolic 
partial differential equation:
\begin{equation}\label{e.FKPPclassical.0}
\frac{\partial u}{\partial t}=\frac{1}{2}\frac{\partial^2u}{\partial x^2}+\beta(u^2-u)
\end{equation}
with $u\in C^{1,2}(\rpn\times\R,[0,1])$. This equation, which first arose in the context of a genetics model for the spread of an advantageous gene through a population, was originally introduced by Fisher \cite{Fis30,Fis37} as well as by Kolmogorov, Petrovskii and Piscounov \cite{KPP37}. Since then it has attracted much attention by analysts and probabilists alike. 
In fact, several authors showed that this equation is closely related to dyadic branching Brownian motions, e.g. \cite{McK75} (see also \cite{McK76}), thus establishing a link of this analytic problem to probability theory. In this probabilistic interpretation the term ``$\frac{1}{2}\frac{\partial^2u}{\partial x^2}$'' corresponds to the motion of the underlying Brownian motion, the ``$\beta$'' is the rate at which the particles split and the term ``$u^2-u$'' results from the binary branching, where two particles replace one particle at each branching time.

A solution $u$ of equation (\ref{e.FKPPclassical.0}) can be interpreted in different ways. The classical work concerning this partial differential equation, such as \cite{Fis30}, \cite{Fis37} and \cite{KPP37}, describes the wave of advance of advantageous genes. More precisely, there are two types of individuals (or genes) in a population and $u(t,x)$
measures the frequency or concentration of the advantageous type at the time-space point $(t,x)$. 
In the  setting regarding the abovementioned probabilistic interpretation, that links
(\ref{e.FKPPclassical.0}) with a dyadic branching Brownian motion, let $u(t,x)$ be the
probability that at time $t$ the largest particle of the branching Brownian motion has a value less than $x$. Then $u$ satisfies
equation (\ref{e.FKPPclassical.0}), see (7) in \cite{McK75}. That is to say, in \cite{Fis30}, \cite{Fis37} and \cite{KPP37} the FKPP equation~(\ref{e.FKPPclassical.0}) describes the bulk of a population, in \cite{McK75} it describes the most advanced particle of a branching Brownian motion.

The classical FKPP travelling waves are solutions of (\ref{e.FKPPclassical.0}) of the form 
\[
u(t,x)=f(x-c t)
\]
for some $f\in C^2(\R,[0,1])$ and some constant $c\in\R$. This leads to the so-called FKPP travelling wave equation 
with wave speed $c\in\R$,
\begin{align*}\label{e.FKPPclassical}
\frac{1}{2}f''+cf'+\beta(f^2-f) &= 0\notag
\\[0.5ex]
\lim_{x\to-\infty}f(x) &= 0
\\[0.5ex]
\lim_{x\to\infty}f(x) &= 1,\notag
\end{align*}
where $\beta>0$. This travelling wave boundary value problem was studied by various authors, using both analytic as well as probabilistic techniques, and it is known that is has a unique (up to additive translation) solution $f\in C^2(\R,[0,1])$ if $|c|\ge\sqrt{2\beta}$. In the opposite case that $0\le|c|<\sqrt{2\beta}$ there is no travelling wave solution. Regarding probabilistic approaches to the classical FKPP travelling wave equation  we also refer for instance to the work of Bramson \cite{Bra78,Bra83}, Chauvin and Rouault \cite{CR88,CR90}, Uchiyama \cite{Uch77,Uch78} as well as \cite{Har99}, \cite{101} and \cite{Nev87}.

Interesting with regard to our work is that the above boundary value problem was extended to continuous-time branching random walks, cf. \cite{Kyp99}, and  to conservative homogeneous fragmentation processes, see \cite{104}. In the context of such fragmentation processes the corresponding partial integro-differential equation, referred to as FKPP equation, is given by
\begin{equation}\label{e.FKPP.0.a}
\frac{\partial u}{\partial t}(t,x)=\int_{\mathcal S}\left(\prod_{n\in\N}u(t,x+\ln(s_n))-u(t,x)\right)\nu(\dd{\bf s}) 
\end{equation}
for certain $u:\rpn\times\R\to[0,1]$. Note that (\ref{e.FKPP.0.a}) looks quite different compared to the classical FKPP equation (\ref{e.FKPPclassical.0}). This difference results from the fact that fragmentation processes have no spatial motion except at jump times and from the more complicated jump structure of fragmentations in comparison with dyadic branching Brownian motions. However, despite the difference of the above equation compared to the classical FKPP equation, the name {\it FKPP equation} for (\ref{e.FKPP.0.a}) stems from the similarity in terms of the probabilistic interpretation these two equations have.
Of particular interest to us are the  FKPP travelling waves to (\ref{e.FKPP.0.a}) with wave speed $c\in\R$, i.e. solutions of (\ref{e.FKPP.0.a}) which are of the form $u(t,x)=f(x-ct)$ for all $t\in\rpn$ and $x\in\R$. These travelling wave solutions are functions $f\in C^1(\R,[0,1])$ that satisfy the following FKPP travelling wave equation 
\[
cf'(x)+\int_{\mathcal S}\left(\prod_{n\in\N}f(x+\ln(s_n)-f(x)\right)\nu(\dd{\bf s})=0
\]
for all $x\in\R$ with boundary conditions
\[
\lim_{x\to-\infty}f(x)=0\qquad\text{and}\qquad\lim_{x\to\infty}f(x)=1.
\]

\begin{rem}\label{r.BHK09}
For every $p\in(\underline{p},\bar p]$ let $\mathcal T_2(p)$ denote the space of monotonically increasing functions $f\in C^1(\R,[0,1])$ satisfying the boundary conditions  $\lim_{x\to-\infty}f(x)=0$ as well as $\lim_{x\to\infty}f(x)=1$ and such that the mapping $x\mapsto e^{(1+p)x}(1-f(x))$ is monotonically increasing. In Theorem~1 of \cite{104} Berestycki et. al. showed that for $p\in(\underline{p},\bar p]$
there exists a unique (up to additive translation) FKPP travelling wave solution in $\mathcal T_2(p)$ with wave speed $c_p$, cf. (\ref{c_p.2}). According to Lemma~1 of \cite{89} the mapping $p\mapsto\nicefrac{\Phi(p)}{(1+p)}=c_p$ is monotonically increasing on $(\underline p,\bar p]$ and thus it follows in view of Theorem~3 (ii) in \cite{104} that  $c_{\bar p}$ is the maximal wave speed for two-sided travelling waves. 
\end{rem}

In this paper we are interested in the one-sided counterpart of the abovementioned FKPP equation. In the classical setting the one-sided FKPP equation 
is the following partial differential equation
\[
\frac{\partial u}{\partial t}=\frac{1}{2}\frac{\partial^2u}{\partial x^2}+\beta(u^2-u)
\]
on $\rp\times\rp$ with $u\in C^{1,2}(\rpn\times\rpn,[0,1])$. Observe that this equation is the analogue of (\ref{e.FKPPclassical.0}) for functions defined on $\rpn\times\rpn$. The corresponding one-sided FKPP travelling wave equation 
with wave speed $c\in\R$ is given by the differential equation
\begin{equation}\label{e.FKPPclassical.2}
\frac{1}{2}f''+cf'+\beta(f^2-f) = 0
\end{equation}
on $\rp$ for $f\in C^2(\rpn,[0,1])$ satisfying the boundary conditions
\begin{equation}\label{e.FKPPclassical.2.0a}
\lim_{x\to0}f(x) = 1\qquad\text{as well as}\qquad\lim_{x\to\infty}f(x) = 0.
\end{equation}
By considering killed branching Brownian motion with drift, killed upon hitting the origin, Harris et. al. proved in \cite{HHK06} that solutions of the one-sided FKPP travelling wave boundary value problem (\ref{e.FKPPclassical.2}) and (\ref{e.FKPPclassical.2.0a}) exist and are unique (up to additive translation) for all  $c\in(-\sqrt{2\beta},\infty)$ and there is no such travelling wave solution for $c\in(-\infty,-\sqrt{2\beta}]$. Notice that the one-sided travelling wave solutions for negative $c$ are precisely those wave speeds for which there does not exist a two-sided travelling wave. With regard to the one-sided FKPP travelling wave equation in the classical setting we refer also to  \cite{Wat65}, concerning existence of a solution, as well as \cite{Pin95} for existence and uniqueness of a solution of (\ref{e.FKPPclassical.2}) and (\ref{e.FKPPclassical.2.0a}) obtained by means of analytic techniques.

\section{The finite activity case}\label{s.fa}
In this section we prove existence of a one-sided FKPP travelling wave in the situation of a finite dislocation measure $\nu$, sometimes referred to as the {\it finite activity case}. In this respect note that a homogeneous fragmentation process with finite $\nu$ may still have infinitely many jumps in any finite time interval after the first jump, because infinitely many blocks may be present at any such time and each block fragments with the same rate. However, in this setting of a finite dislocation measure every block $B_n$, $n\in\N$, has only finitely many jumps up to any $t\in\rpn$. In particular, this implies that the fragmentation process with killing and with finite $\nu$ has finite activity in bounded time intervals, since at any time there are only finitely many blocks alive. Therefore, in this finite activity situation it is possible to consider the time $\tau_k:\Omega\to\rp\cup \{\infty\}$, $k\in\N$, of the $k$-th jump of the killed process $\Pi^x$, a fact we shall make use of below.

An approach to solve the classical one-sided FKPP equation with boundary condition $u(0,x)=g(x)$ for some suitable function $g:\rpn\to[0,1]$ is to show that the function $u:\rpn\times\rpn\to[0,1]$, given by
\[
u(t,x)=\mathbb E\left(\prod_{n\in\N}g(x+Y_n(t))\right)
\]
for all $t,x\in\rpn$, is a solution of the considered boundary value problem, where the $Y_n(t)$
are the positions of the particles at time $t$ in a dyadic branching Brownian motion. 

In this section we show that for fragmentations with a finite dislocation measure $\nu$ a similar approach as above works for the initial value problem (\ref{e.FKPP.0}). 
More precisely, we prove that for certain functions $g:\R_\infty\to[0,1]$ the function $u:\rpn\times\R_\infty\to[0,1]$, defined by 
\begin{equation}\label{e.classicalFKPP}
\forall\,x\in[ct,\infty):\,u(t,x)=\mathbb E\left(\prod_{n\in\mathcal N^{x-ct}_t}g(x+\ln(|B_n(t)|))\right)\qquad\text{and}\qquad u(t,\cdot)|_{[-\infty,ct)}\equiv1
\end{equation}
for all $t\in\rpn$,  solves equation~(\ref{e.FKPP.0}) with boundary condition $u(0,\cdot)=g$. 

\begin{proposition}\label{p.t.mainresult.1.1}
Assume that $\nu(\mathcal S)<\infty$ and let $c>0$. Then every function $u:\rpn\times\R_\infty\to[0,1]$ defined by (\ref{e.classicalFKPP}), for some function $g:\R_\infty\to[0,1]$ with $g|_{\rpn}\in C^0(\rpn)$ and $g|_{[-\infty,0)}\equiv1$,  satisfies  the boundary condition
\begin{equation}\label{e.bc}
u(0,\cdot)=g
\end{equation} and solves equation~(\ref{e.FKPP.0}) for any $x\in\rpn$. In particular, any such function $u$ of the form $u(t,x)=f(x-ct)$ for some $f:\R_\infty\to[0,1]$, with $f|_{[-\infty,0)}\equiv1$, and all $t\in\rpn$ and $x\in\R_\infty$ solves (\ref{e.FKPP.1}).
\end{proposition}

The following corollary of Proposition~\ref{p.t.mainresult.1.1} provides a  short proof that the extinction probability $\varphi$ of the  fragmentation process with killing solves equation~(\ref{e.FKPP.1}) in the special case of a finite dislocation measure.

\begin{corollary}\label{c.p.t.mainresult.1.1.1}
Assume that $\nu(\mathcal S)<\infty$ and let $c>c_{\bar p}$. Then $\varphi$ is an FKPP travelling wave with wave speed $c$.
\end{corollary}

\begin{pf}
Let us first show that $\varphi$ solves (\ref{e.FKPP.1}). For this purpose, observe that the fragmentation property, in conjunction with the tower property of conditional expectations, yields that
\begin{align*}
\varphi(x-ct) &= \mathbb E\left(\prod_{n\in\mathcal N^{x-ct}_t}\left.\mathbb P\left(\zeta^{x-ct+ct+y}<\infty\right)\right|_{y=\ln(|B_n(t)|)}\right)
\\[0.5ex]
&= \mathbb E\left(\prod_{n\in\mathcal N^{x-ct}_t}\varphi(x+\ln(|B_n(t)|))\right),
\end{align*}
and thus $u:\rpn\times\R_\infty\to[0,1]$, given by $u(t,x):=\varphi(x-ct)$, satisfies (\ref{e.classicalFKPP}) with $g=\varphi$. Hence, according to Proposition~\ref{p.t.mainresult.1.1} the function $\varphi$ solves (\ref{e.FKPP.1}). Since $c>c_{\bar p}$, it follows from Theorem~10 in \cite{KK12} that $\varphi$ also satisfies the boundary condition~(\ref{e.FKPP.2}), which completes the proof.
\end{pf}

The major part of this paper, cf. Theorem~\ref{t.mainresult.1}, is concerned with the proof that the conclusion of Corollary~\ref{c.p.t.mainresult.1.1.1} holds true also in the general case of an infinite dislocation measure.

{\bf Proof of Proposition~\ref{p.t.mainresult.1.1}}\quad The proof is based on a decomposition according to the first and second jump times of a killed fragmentation. Treating these parts separately we obtain the desired expression for the right derivative of $u(\cdot,x)$, $x\in\rpn$. 

Let $g:\R_\infty\to[0,1]$ be some  function that satisfies $g|_{\rpn}\in C^0(\rpn)$ and $g|_{[-\infty,0)}\equiv1$. Further, consider the function $u:\rpn\times\R_\infty\to[0,1]$ defined by (\ref{e.classicalFKPP}) and fix some $x\in\rpn$ as well as $t\in\mathcal C_{u(\cdot,x)}$. In the light of the c\`adl\`ag paths of $\Pi$ and the DCT, note first that $u$ satisfies the boundary condition (\ref{e.bc}), since $|B_1(0)|=1$ and $|B_n(0)|=0$, i.e. $g(x+\ln(|B_n(0)|))=1$, for all $n\in\N\setminus\{1\}$.  In order to prove that $u$ solves (\ref{e.FKPP.0}) Lebesgue-a.e.,  let $(t_i)_{i\in\mathcal I^x}$ be the jump times of $\Pi^x$ and in view of the finiteness of the dislocation measure and $N^x_t\le e^{x+ct}$, cf. (\ref{e.N_x}), we assume without loss of generality that $\mathcal I^x=\N$ and that $0<t_i< t_j$ for any $i,j\in\N$ with $i<j$.   Since $t_1$ is exponentially distributed with parameter $\mu(\mathcal P)$, we have
\[
\lim_{h\downarrow0}\frac{\mathbb P\left(t_1\le h\right)}{h}\notag
= \lim_{h\downarrow0}\frac{1-e^{-h\mu(\mathcal P)}}{h}
= \mu(\mathcal P)
\]
and deduce by resorting to the strong fragmentation property of $\Pi$ that
\begin{equation}\label{second_jump}
\begin{aligned}
\lim_{h\downarrow0}\frac{\mathbb P(t_2\le h)}{h}
&\le \lim_{h\downarrow0}\frac{\mathbb P(t_1\le h)}{h}\lim_{h\downarrow0}\mathbb E\left(\left.\mathbb P\left(\mathfrak{e}_{\mu(\mathcal P)e^{x+c t}}\le h\right)\right|_{t=t_1}\right)
\\[0.5ex]
&= \mu(\mathcal P)\mathbb E\left(\lim_{h\downarrow0}\left(1-e^{-h\mu(\mathcal P)e^{x+ct_1}}\right)\right)
\\[0.5ex]
&=0,
\end{aligned}
\end{equation}
where  $\mathfrak{e}_{\mu(\mathcal P)e^{x+c t_1}}$ denotes a random variable that is exponentially distributed with parameter $\mu(\mathcal P)e^{x+c t_1}$.
Consequently,
\begin{equation}\label{e.p.t.mainresult.1.1.3}
\lim_{h\downarrow0}\frac{\mathbb P(t_1\le h<t_2)}{h}
= \lim_{h\downarrow0}\frac{\mathbb P(t_1\le h)}{h}-\lim_{h\downarrow0}\frac{\mathbb P(t_2\le h)}{h}
= \mu(\mathcal P).
\end{equation}
By means of the strong Markov property, the fact that the distrubution of $\pi(t_1)$ is given by $\nicefrac{\mu(\cdot)}{\mu(\mathcal P)}$ and the independence between $\pi(t_1)$ and the random vector $
\left(
\begin{array}{c}
t_1
\\
t_2
\end{array}
\right)$, see Proposition~2 in Section 0.5 of \cite{Ber96},  we have 
\begin{align*}
&  \mathbb E\left(\prod_{n\in\mathcal N^{x-c(t+h)}_{t+h}}g(x+\ln(|B_n(t+h)|))\mathds1_{\left\{t_1\le h<t_2\right\}}\right)
\\[0.75ex]
&= \mathbb E\Bigg(\prod_{n\in\N}\Bigg(\mathds1_{\{x-c(t+h)+ct_1<-\ln(|\pi_n(t_1)|)\}}+\mathds1_{\{x-c(t+h)+ct_1\ge-\ln(|\pi_n(t_1)|)\}}
\\[0.5ex]
&\qquad \cdot\mathbb E\Bigg(\prod_{k\in\mathcal N^{x-c(t+h)+ct_1+\ln(|\pi_n(t_1)|)}_t}g\left(x+\ln(|\pi_n(t_1)|)+\ln\left(\left|B^{(n)}_k(t)\right|\right)\right)\mathds1_{\left\{t_1\le h<t_2\right\}}\Bigg|\mathscr F_{t_1}\Bigg)\Bigg)\Bigg)
\\[0.75ex]
&= \mathbb E\Bigg(\prod_{n\in\N}\Bigg(\mathds1_{\{x-c(t+h-t_1)<-\ln(|\pi_n(t_1)|)\}}+\mathds1_{\{x-c(t+h-t_1)\ge-\ln(|\pi_n(t_1)|)\}}
\\[0.5ex]
&\qquad \cdot\mathbb E\Bigg(\mathds1_{\left\{t_1\le h<t_2\right\}}\mathbb E\Bigg(\prod_{k\in\mathcal N^{x-c(t+h-t_1)+\ln(|\pi_n(t_1)|)}_t}g\left(x+\ln(|\pi_n(t_1)|)+\ln\left(\left|B^{(n)}_k(t)\right|\right)\right)\Bigg|\mathscr F_h\Bigg)\Bigg|\mathscr F_{t_1}\Bigg)\Bigg)\Bigg)
\\[0.75ex]
&= \mathbb E\Bigg(\mathds1_{\left\{t_1\le h<t_2\right\}}\prod_{n\in\N}\Bigg(\mathds1_{\{x-c(t+h-t_1)<-\ln(|\pi_n(t_1)|)\}}+\mathds1_{\{x-c(t+h-t_1)\ge-\ln(|\pi_n(t_1)|)\}}
\\[0.5ex]
&\qquad \cdot\mathbb E\Bigg(\prod_{k\in\mathcal N^{x+\ln(u)-ct}_t}g\left(x+\ln\left(u\right)+\ln(|B_k(t)|)\right)\Bigg)\Bigg|_{u=|\pi_n(t_1)|}\Bigg)\Bigg)
\\[0.75ex]
&=  \mathbb P(t_1\le h<t_2)\mathbb E\left(\prod_{n\in\N}u_h(t,x+\ln(|\pi_n(t_1)|))\right)
\\[0.75ex]
&=  \mathbb P(t_1\le h<t_2)\int_\mathcal P\prod_{n\in\N}u_h(t,x+\ln(|\pi_n|))\frac{\mu(\dd\pi)}{\mu(\mathcal P)}\, ,
\end{align*}
where
\[
u_h(t,\cdot)|_{[c(t+h-t_1),\infty)}:= u|_{[c(t+h-t_1),\infty)}\qquad\text{as well as}\qquad u_h(t,\cdot)|_{[-\infty,c(t+h-t_1))}:\equiv1.
\]
Therefore, (\ref{e.p.t.mainresult.1.1.3}) yields that
\begin{align}\label{e.p.t.mainresult.1.1.3_b}
& \lim_{h\downarrow0}\mathbb E\left(\frac{1}{h}\prod_{n\in\mathcal N^{x-c(t+h)}_{t+h}}g(x+\ln(|B_n(t+h)|))\mathds1_{\left\{t_1\le h<t_2\right\}}\right) \notag
\\[0.75ex]
&= \int_\mathcal P\prod_{n\in\N}\lim_{h\downarrow0}u_h(t,x+\ln(|\pi_n|))\mu(\dd\pi)\frac{1}{\mu(\mathcal P)}\lim_{h\downarrow0}\frac{\mathbb P(t_1\le h<t_2)}{h}
\\[0.75ex]
&= \int_\mathcal P\prod_{n\in\N}u(t,x+\ln(|\pi_n|))\mu(\dd\pi).\notag
\end{align}
Moreover, 
\begin{align*}
&\mathbb E\left(\prod_{n\in\mathcal N^{x-ct}_t}g(x+\ln(|B_n(t)|))\mathds1_{\left\{t^{(n)}_1\le h<t^{(n)}_2\right\}}\right)
\\[0.5ex]
&= \mathbb E\left(\prod_{n\in\mathcal N^{x-ct}_t}g(x+\ln(|B_n(t)|))\mathbb P\left(\left.t^{(n)}_1\le h<t^{(n)}_2\right|\mathscr F_t\right)\right)
\\[0.5ex]
&=u(t,x)\mathbb P(t_1\le h<t_2)
\\[0.5ex]
&=\int_\mathcal Pu(t,x)\mu(\dd\pi)\frac{1}{\mu(\mathcal P)}\mathbb P(t_1\le h<t_2)
\end{align*}
holds for all $h>0$, where conditionally on $\mathscr F_t$ the $t^{(n)}_1$ and $t^{(n)}_2$ are independent copies of $t_1 $ and $t_2$, respectively. 
Hence,
\begin{align}\label{e.p.t.mainresult.1.1.3_c}
\lim_{h\downarrow0}\mathbb E\left(\frac{1}{h}\prod_{n\in\mathcal N^{x-ct}_t}g(x+\ln(|B_n(t)|))\mathds1_{\left\{t^{(n)}_1\le h<t^{(n)}_2\right\}}\right) \notag
&= \int_\mathcal Pu(t,x)\mu(\dd\pi)\frac{1}{\mu(\mathcal P)}\lim_{h\downarrow0}\frac{\mathbb P(t_1\le h<t_2)}{h}
\\[0.5ex]
&= \int_\mathcal Pu(t,x)\mu(\dd\pi).
\end{align}
Since $|B_1(h)|=1$ and $\mathcal N^{x-c(t+h)}_h=\{1\}$ on $\{t_1> h\}$, we deduce with $t^{(n)}_1$ being defined as above that
\begin{equation}\label{no_jump}
\begin{aligned}
& \mathbb E\left(\prod_{n\in\mathcal N^{x-c(t+h)}_{t+h}}g(x+\ln(|B_n(t+h)|))\mathds1_{\left\{t_1> h\right\}}\right)
\\[0.75ex]
&= \mathbb E\left(\left.\mathbb E\left.\left(\prod_{n\in\tilde{\mathcal N}^{x-c(t+h)+ch}_t}g(x+\ln(\gamma)+\ln(|B^{(n)}(t)|))\right|\mathscr F_h\right)\right|_{\gamma=|B_1(h)|}\mathds1_{\left\{t_1> h\right\}}\right)
\\[0.75ex]
&= \mathbb E\left(\mathbb E\left(\prod_{n\in\mathcal N^{x-ct}_t}g(x+\ln(|B_n(t)|))\right)\mathds1_{\left\{t_1> h\right\}}\right)
\\[0.75ex]
&= \mathbb E\left(\prod_{n\in\mathcal N^{x-ct}_t}g(x+\ln(|B_n(t)|))\right)\mathbb P\left(t_1> h\right)
\\[0.75ex]
&= \mathbb E\left(\prod_{n\in\mathcal N^{x-ct}_t}g(x+\ln(|B_n(t)|))\mathbb P\left(t_1> h\right)\right)
\\[0.75ex]
&= \mathbb E\left(\prod_{n\in\mathcal N^{x-ct}_t}g(x+\ln(|B_n(t)|))\mathbb P\left(\left.t^{(n)}_1> h\right|\mathscr F_t\right)\right)
\\[0.75ex]
&=\mathbb E\left(\prod_{n\in\mathcal N^{x-ct}_t}g(x+\ln(|B_n(t)|))\mathds1_{\left\{t^{(n)}_1> h\right\}}\right)
\end{aligned}
\end{equation}
holds for each $h>0$, where conditionally on $\mathscr F_h$ the $\tilde{\mathcal N}^{(\cdot)}_t$ and $B^{(n)}$ are independent copies of $\mathcal N^{(\cdot)}_t$ and $B_n$, respectively. Furthermore, note that (\ref{second_jump}) results in
\begin{equation}\label{second_jump_b}
\lim_{h\downarrow0}\frac{\mathbb P\left(t^{(n)}_2\le h\right)}{h}=\lim_{h\downarrow0}\frac{\mathbb E\left(\mathbb P\left(\left.t^{(n)}_2\le h\right|\mathscr F_t\right)\right)}{h}=\lim_{h\downarrow0}\frac{\mathbb P\left(t_2\le h\right)}{h}=0,
\end{equation}
where $t^{(n)}_2$ is defined as above.
Bearing in mind that $|g|\le1$ it follows from the DCT in conjunction with (\ref{second_jump}) and (\ref{second_jump_b}), respectively, that
\begin{align*}
\lim_{h\downarrow 0}\mathbb E\left(\prod_{n\in\mathcal N^{x-c(t+h)}_{t+h}}g(x+\ln(|B_n(t+h)|))\mathds1_{\left\{t_2\le h\right\}}\right)
&=0
\\[0.5ex]
&=\lim_{h\downarrow 0}\mathbb E\left(\prod_{n\in\mathcal N^{x-ct}_t}g(x+\ln(|B_n(t)|))\mathds1_{\left\{t^{(n)}_2\le h\right\}}\right).
\end{align*}
Since
\begin{align*}
u(t+h,x) &= \mathbb E\left(\prod_{n\in\mathcal N^{x-c(t+h)}_{t+h}}g(x+\ln(|B_n(t+h)|))\mathds1_{\left\{t_1>h\right\}}\right)
\\[0.75ex]
&\qquad +\mathbb E\left(\prod_{n\in\mathcal N^{x-c(t+h)}_{t+h}}g(x+\ln(|B_n(t+h)|))\mathds1_{\left\{t_2\le h\right\}}\right)
\\[0.75ex]
&\qquad +\mathbb E\left(\prod_{n\in\mathcal N^{x-c(t+h)}_{t+h}}g(x+\ln(|B_n(t+h)|))\mathds1_{\left\{t_1\le h<t_2\right\}}\right)
\end{align*}
and
\begin{align*}
u(t,x) &= \mathbb E\left(\prod_{n\in\mathcal N^{x-ct}_t}g(x+\ln(|B_n(t)|))\mathds1_{\left\{t^{(n)}_1> h\right\}}\right)+ \mathbb E\left(\prod_{n\in\mathcal N^{x-ct}_t}g(x+\ln(|B_n(t)|))\mathds1_{\left\{t^{(n)}_2\le h\right\}}\right)
\\[0.5ex]
&\qquad +\mathbb E\left(\prod_{n\in\mathcal N^{x-ct}_t}g(x+\ln(|B_n(t)|))\mathds1_{\left\{t^{(n)}_1\le h<t^{(n)}_2\right\}}\right)
\end{align*}
hold for every $h>0$, it thus follows from (\ref{e.p.t.mainresult.1.1.3_b}), (\ref{e.p.t.mainresult.1.1.3_c}) and (\ref{no_jump}) that
\[
\lim_{h\downarrow 0}\frac{u(t+h,x)-u(t,x)}{h}
= \int_\mathcal P\left(\prod_{n\in\N}u(t,x+\ln(|\pi_n|))-u(t,x)\right)\mu(\dd\pi),
\]
which completes the proof, since $t\in\mathcal C_{u(\cdot,x)}$.
\hfill$\square$

\section{Sufficiency criterion for the existence of travelling waves}\label{s.cspm}

The goal of this section is to provide the proof of Theorem~\ref{p.t.1.2}. 
A first approach to try proving Theorem~\ref{p.t.1.2} might be to pursue a line of argument along the lines of the proof of Theorem~1 in \cite{104}. But that proof relies on $f$ being continuously differentiable and in our situation we cannot use any differentiability assumption. In fact, even if we knew that $f$ is differentiable with a bounded derivative $f'$, we would at least need that the set of discontinuities of $f'$ is a Lebesgue null set. However, in general the set of such discontinuities may have positive Lebesgue measure, cf. Example~3.5 in \cite{119}. 

Let us start with the following auxiliary result.

\begin{lemma}\label{l.MVT}
Let $f\in\mathcal T$ as well as $a,b\in\rp$. Then we have
\[
f(a)-f(b)\le(b-a)\sup_{x\in(a,b)}|f'_+(x)|.
\]
\end{lemma}

\begin{pf}
Define a function $\phi:[a,b]\to\R$ by
\[
\phi(x):=f(x)-\frac{f(b)-f(a)}{b-a}(x-a)
\]
for all $x\in[a,b]$. 
Let us first show that there exists some $x_0\in(a,b)$ such that
\begin{equation}\label{e.MVT.-1}
\limsup_{h\downarrow0}\frac{\phi(x_0+h)-\phi(x_0)}{h}\le0.
\end{equation}
To this end, assume
\begin{equation}\label{e.MVT.0}
\phi'_+(x):=\limsup_{h\downarrow0}\frac{\phi(x+h)-\phi(x)}{h}>0
\end{equation}
for each $x\in(a,b)$. Then for every $x\in(a,b)$ there exists some $\epsilon_x>0$ such that for every $\epsilon\in(0,\epsilon_x]$ we have
\begin{equation}\label{e.MVT.1}
\frac{\phi(x+h)-\phi(x)}{h}>0
\end{equation}
for some $h\in(0,\epsilon)$.
We now show that this implies that $\phi$ is nondecreasing on $(a,b)$. For this purpose, consider $c,d\in[a,b]$ and assume
\begin{equation}\label{e.MVT.2}
\max_{x\in[c,d]}\phi(x)\ne\phi(d),
\end{equation}
where the existence of this maximum follows from the continuity of $\phi$, which in turn follows from $f\in\mathcal T$ being continuous. Then there exists some $x_0\in[c,d)$ such that
\[
\max_{x\in[c,d]}\phi(x)=\phi(x_0).
\]
However, this implies that $\phi(x_0)\ge\phi(x)$ for all $x\in(x_0,(x_0+\epsilon_{x_0})\land d)$, which contradicts (\ref{e.MVT.1}). Hence, (\ref{e.MVT.2}) cannot be true and consequently we infer that
\begin{equation}\label{e.MVT.3}
\max_{x\in[c,d]}\phi(x)=\phi(d)
\end{equation}
for all $c,d\in[a,b]$ under assumption~(\ref{e.MVT.0}). Note that $\phi$ not being nondecreasing on $[a,b]$ would entail that there exist $c,d\in[a,b]$, with $c<d$, such that $\phi(c)>\phi(d)$, which contradicts (\ref{e.MVT.3}). Therefore, we conclude that $\phi$ is nondecreasing and nonconstant on $[a,b]$ if (\ref{e.MVT.0}) holds. This, however, contradicts the fact that
\[
\phi(a)=f(a)=\phi(b).
\]
We thus deduce that (\ref{e.MVT.0}) cannot hold and hence there exists some $x_0\in(a,b)$ such that (\ref{e.MVT.-1}) holds.

With $x_0\in(a,b)$ given by (\ref{e.MVT.-1}) we  obtain
\[
0\ge\phi'_+(x_0)=f'_+(x_0)-\frac{f(b)-f(a)}{b-a},
\]
which results in
\[
0\le \sup_{x\in(a,b)}|f'_+(x)|-\frac{f(a)-f(b)}{b-a}
\]
and thus
\[
f(a)-f(b)\le (b-a)\sup_{x\in(a,b)}|f'_+(x)|.
\]
\end{pf}

We proceed by establishing two  auxiliary results, which in spirit are analogues of respective results in \cite{104}. Afterwards we provide a lemma giving conditions under which only the block containing $1$ is alive in the  fragmentation process with killing. Finally, having all these auxiliary results at hand, we finish this section with the proof of Theorem~\ref{p.t.1.2}.

Observe first that a straightforward argument by induction yields that
\begin{equation}\label{e.estimate}
\left|\prod_{n\in\N}a_n-\prod_{n\in\N}b_n\right|\le\sum_{n\in\N}|a_n-b_n|
\end{equation}
holds for all sequences $(a_n)_{n\in\N},(b_n)_{n\in\N}\in[0,1]^\N$. The following lemma, whose proof is based on (\ref{e.estimate}), shows in particular that $f\in\mathcal D_L$ and, moreover, that $Lf$ is bounded on compact sets for any $f\in\mathcal T$.

\begin{lemma}\label{l.p.t.1.2.4.0}
Let $f\in\mathcal T$ and let $a,b\in\rp$. Then
\[
\int_{\mathcal P}\sup_{x\in[a,b]}\left|\prod_{n\in\N}f(x+\ln(|\pi_n|))-f(x)\right|\mu(\dd\pi)<\infty.
\] 
\end{lemma}

\begin{pf}
By means of (\ref{e.estimate}) we have
\begin{align}\label{e.l.p.t.1.2.3.1}
& \int_{\mathcal P}\sup_{x\in[a,b]}\left|\prod_{n\in\N}f(x+\ln(|\pi_n|))-f(x)\right|\mu(\dd\pi)
\\[0.5ex]
&\le \int_{\mathcal P}\sup_{x\in[a,b]}\left|f(x+\ln(|\pi|^\downarrow_1))-f(x)\right|\mu(\dd\pi)+\int_{\mathcal P}\sum_{n\in\N\setminus\{1\}}\sup_{x\in[a,b]}\left|f(x+\ln(|\pi|^\downarrow_n))-1\right|\mu(\dd\pi).\notag
\end{align}

Since 
\[
\frac{\dd}{\dd x}[\ln(x)+2(1-x)]=\frac{1}{x}-2
\]
and $\ln(1)+2(1-1)=0$, we deduce that 
\[
-\ln(x)\le2(1-x)
\]
holds for all $x\in[\nicefrac{1}{2},1]$. Therefore, for every $\epsilon\in(0,\nicefrac{1}{2}]$ we have
\begin{equation}\label{e.l.p.t.1.2.3.1a}
-\ln(|\pi|^\downarrow_1)\le2(1-|\pi|^\downarrow_1)
\end{equation} 
for all $\pi\in\mathcal P$ with $1-|\pi|^\downarrow_1\le\epsilon$. Moreover, by means of Lemma~\ref{l.MVT} we have for any $x\in\rp$ and $\pi\in\mathcal P$ with $|\pi|^\downarrow_1>e^{-x}$ the estimate
\begin{equation}\label{e.l.p.t.1.2.3.1b2}
\left|f\left(x+\ln(|\pi|^\downarrow_1)\right)-f(x)\right|\le-\ln(|\pi|^\downarrow_1)\sup_{y\in\left(x+\ln(|\pi|^\downarrow_1),\,x\right)}|f'_+(y)|.
\end{equation}
Furthermore, for every $\gamma\in(0,a)$ define
\[
A_{a,\gamma}:=\left\{\pi\in\mathcal P:a+\ln(|\pi|^\downarrow_1)\in[0,\gamma)\right\}=\left\{\pi\in\mathcal P:|\pi|^\downarrow_1\in[e^{-a},e^{\gamma-a})\right\}
\]
and observe that in view of $\gamma-a<0$ and (\ref{e.levymeasure}) we have
\[
\mu\left(A_{a,\gamma}\right)\le\mu\left(\{\pi\in\mathcal P:|\pi|^\downarrow_1<e^{\gamma-a}\}\right)<\infty.
\]
Hence, resorting to (\ref{e.levymeasure}), (\ref{finite_derivative}) and (\ref{e.l.p.t.1.2.3.1a}) as well as (\ref{e.l.p.t.1.2.3.1b2}) we conclude in the light of (3) in \cite{HKK10} and $f(x)\in[0,1]$ for every $x>0$ that
\begin{align*}
&\int_{\mathcal P}\sup_{x\in[a,b]}\left|f(x+\ln(|\pi|^\downarrow_1))-f(x)\right|\mu(\dd\pi)\notag 
\\[0.5ex]
&\le \int_{\{\pi\in\mathcal P:1-|\pi|^\downarrow_1>\epsilon\}\cup A_{a,\gamma}}\sup_{x\in[a,b]}\left|f(x+\ln(|\pi|^\downarrow_1))-f(x)\right|\mu(\dd\pi)\notag
\\[0.5ex]
&\qquad +\int_{\{\pi\in\mathcal P:1-|\pi|^\downarrow_1\le\epsilon\}\setminus A_{a,\gamma}}\sup_{x\in[a,b]}\left|f(x+\ln(|\pi|^\downarrow_1))-f(x)\right|\mu(\dd\pi)
\\[0.5ex]
&\le \mu\left(\{\pi\in\mathcal P:1-|\pi|^\downarrow_1>\epsilon\}\cup A_{a,\gamma}\right)+\int_{\{\pi\in\mathcal P:1-|\pi|^\downarrow_1\le\epsilon\}\setminus A_{a,\gamma}}-\ln(|\pi|^\downarrow_1)\sup_{y\in(a+\ln(|\pi|^\downarrow_1),b)}|f_+'(y)|\,\mu(\dd\pi)\notag
\\[0.5ex]
&\le \mu\left(\{\pi\in\mathcal P:|\pi|^\downarrow_1<1-\epsilon\}\right)+\mu\left(A_{a,\gamma}\right)+2\sup_{y\in[\gamma,b)}|f_+'(y)|\int_{\mathcal P}(1-|\pi|^\downarrow_1)\,\mu(\dd\pi)\notag
\\[0.5ex]
&< \infty\notag
\end{align*} 
for any $\epsilon\in(0,\nicefrac{1}{2}]$, which shows that the first term on the right-hand side of (\ref{e.l.p.t.1.2.3.1}) is finite. In order to deal with the second term on the right-hand side of (\ref{e.l.p.t.1.2.3.1}), note that the monotonicity of $f$ together with $f|_{[-\infty,0)}\equiv1$ and $f|_{[0,\infty)}\in[0,1]$ yields that
\begin{align*}
\int_{\mathcal P}\sum_{n\in\N\setminus\{1\}}\sup_{x\in[a,b]}|1-f(x+\ln(|\pi|^\downarrow_n))|\mu(\dd\pi) 
&\le \int_{\mathcal P}\sum_{n\in\N\setminus\{1\}}|1-f(b+\ln(|\pi|^\downarrow_n))|\mu(\dd\pi)
\\[0.5ex]
&\le \int_{\mathcal P}\sum_{n\in\N\setminus\{1\}}e^{(b+\ln(|\pi|^\downarrow_n))}\mu(\dd\pi)\notag
\\[0.5ex]
&= e^b\int_{\mathcal P}\sum_{n\in\N\setminus\{1\}}|\pi|^\downarrow_n\mu(\dd\pi)
\\[0.5ex]
&< \infty
\end{align*}
for all $x>0$. Observe that the finiteness holds, since
\[
\int_{\mathcal P}\sum_{n\in\N\setminus\{1\}}|\pi|^\downarrow_n\mu(\dd\pi) =\int_{\mathcal P}\left(\left(1-|\pi|^\downarrow_1\right)+\left(\sum_{n\in\N}|\pi|^\downarrow_n-1\right)\right)\mu(\dd\pi)
\le \int_{\mathcal P}(1-|\pi|^\downarrow_1)\mu(\dd\pi)
< \infty.
\]
Consequently, also the second term on the right-hand side of (\ref{e.l.p.t.1.2.3.1}) is finite.
\end{pf}

As already mentioned, the previous lemma implies that $Lf$ exists for each $f\in\mathcal T$. The next lemma goes a step further for that it shows that $Lf$ is continuous for every $f\in\mathcal T$. 

\begin{lemma}\label{l.p.t.1.2.4}
Let $f\in\mathcal T$. Then the function $Lf$ is continuous on $\rp$.
\end{lemma}

\begin{pf}
Fix some $x\in\rp$ and let $(x_k)_{k\in\N}$ be a sequence in $\rp$ with $x_k\to x$ as $k\to\infty$. In addition, fix some $\epsilon\in(0,x)$ and let $k_\epsilon\in\N$ be such that $|x-x_k|\le\epsilon$ for all $k\ge k_\epsilon$.
Observe that 
\begin{align}\label{e.l.p.t.1.2.3.3}
& \int_{\mathcal P}\sup_{k\ge k_\epsilon}\left|\prod_{n\in\N}f(x+\ln(|\pi_n|))-f(x)-\prod_{n\in\N}f(x_k+\ln(|\pi_n|))+f(x_k)\right|\mu(\dd\pi)
\\[0.5ex]
&\le \int_{\mathcal P}\left|\prod_{n\in\N}f(x+\ln(|\pi_n|))-f(x)\right|\mu(\dd\pi)
+\int_{\mathcal P}\sup_{y\in[x-\epsilon,\,x+\epsilon]}\left|\prod_{n\in\N}f(y+\ln(|\pi_n|))-f(y)\right|\mu(\dd\pi).\notag
\end{align}
According to Lemma~\ref{l.p.t.1.2.4.0} both of the integrals on the right-hand side of (\ref{e.l.p.t.1.2.3.3}) are finite. Hence, we can apply the DCT and deduce that
\begin{align*}
&\lim_{k\to\infty}|Lf(x)-Lf(x_k)| 
\\[0.5ex]
&= \int_{\mathcal P}\lim_{k\to\infty}\left|\prod_{n\in\N}f(x+\ln(|\pi_n|))-f(x)-\prod_{n\in\N}f(x_k+\ln(|\pi_n|))+f(x_k)\right|\mu(\dd\pi)
\\[0.5ex]
&= \int_{\mathcal P}\left|\prod_{n\in\N}f(x+\ln(|\pi_n|))-\prod_{n\in\N}\lim_{k\to\infty}f(x_k+\ln(|\pi_n|))-f(x)+\lim_{k\to\infty}f(x_k)\right|\mu(\dd\pi)
\\[0.5ex]
&= 0,
\end{align*}
where the final equality follows from  $f\in\mathcal T$ being continuous on $\rp$. Notice that we can interchange the limit and the product in the penultimate equality, since only finitely many factors of the product differ from 1. Hence, we have proven the continuity of $Lf$ at $x$ and since $x\in\rp$ was chosen arbitrarily, this completes the proof.
\end{pf}

Recall the process $Z^{x,f}$ that we defined in (\ref{e.prod_martingale}) and set
\[
\Delta Z^{x,f}_t:=Z^{x,f}_t-Z^{x,f}_{t-}
\]
for every $t>0$.

We are now in a position to prove  Theorem~\ref{p.t.1.2}.

{\bf Proof of Theorem~\ref{p.t.1.2}}
Throughout the proof let $x\in\mathcal C_f$ and let $(a_n)_{n\in\N}$ be a sequence in $(0,1)$ with $a_n\downarrow 0$ as $n\to\infty$. Moreover, Consider the following stopping time
\begin{equation}\label{d.delta}
\delta:=\inf\left\{t>0:x+ct+\ln\left(|\Pi(t)|^\downarrow_2\right)>0\right\}\land\tau^-_{1,0}\land1.
\end{equation}
Recall from Lemma~\ref{l.p.t.1.2.4.0} that $f\in\mathcal D_L$.

The idea of the proof is to consider  an appropriate decomposition of the limit of $\mathbb E(Z^{x,f}_{\delta\land a_n}-Z^{x,f}_{0})a_n^{-1}$ as $n\to\infty$, which by the martingale property of $Z^{x,f}$ equals $0$. In this spirit the proof  deals with the jumps and drift that contribute to the difference $Z^{x,f}_{\delta\land a_n}-Z^{x,f}_{0}$ separately and eventually combines these considerations in order  to prove the assertion.

Let us first deal with the jumps of $\Pi^x$ that contribute to the difference $Z^{x,f}_{\delta\land a_n}-Z^{x,f}_{0}$.
To this end, we start by pointing out that $\delta>0$ $\mathbb P$-almost surely. Indeed, if $\delta=\tau^-_{1,0}\land1$, then the $\mathbb P$-a.s. positivity of $\delta$ follows, since for $X_n$ the point $0$ is irregular for $(-\infty,0)$. In order to deal with the case $\delta<\tau^-_{1,0}\land1$, note that
\[
x+c\tau+\ln\left(|\Pi(\tau)|^\downarrow_1\right)\ge x+c\tau+\ln\left(|B_1(\tau)|\right)=X^x_1(\tau)\ge x\Longrightarrow|\Pi(\tau)|^\downarrow_1\ge e^{-c\tau}
\]
and
\[
x+c\tau+\ln\left(|\Pi(\tau)|^\downarrow_2\right)\ge0
\Longrightarrow|\Pi(\tau)|^\downarrow_2\ge e^{-(x+c\tau)}
\]
hold on the event $\{\delta<\tau^-_{1,0}\land1\}$ 
for any random time $\tau\le\delta$. Therefore, on this event
we have
\[
e^{-c\tau}\le|\Pi(\tau)|^\downarrow_1\le1-e^{-(x+c\tau)}
\]
for every random time $\tau\le\delta$, which implies that
\[
\delta\ge\tau\ge\frac{1}{c}\ln\left(1+e^{-x}\right)>0
\]
on $\{\delta<\tau^-_{1,0}\land1\}$.

The compensation formula for Poisson point processes yields that
\begin{align*}
& \frac{1}{a_n}\mathbb E\left(\sum_{i\in\mathcal I}\mathds1_{(0,\delta\land a_n]}(t_i)\Delta Z^{x,f}_{t_i}\right)
\\[0.5ex]
&= \frac{1}{a_n}\mathbb E\left(\sum_{i\in\mathcal I}\mathds1_{(0,\delta\land a_n]}(t_i)\prod_{l\in\N}f(X^x_1(t_i-)+\ln(|\pi_l(t_i)|))-f(X^x_1(t_i-))\right)
\\[0.5ex]
&= \frac{1}{a_n}\mathbb E\left(\int_{(0,\delta\land a_n]}\int_\mathcal P\prod_{l\in\N}f(X^x_1(t-)+\ln(|\pi_l|))-f(X^x_1(t-))\mu(\dd\pi)\dd t\right)
\\[0.5ex]
&= \mathbb E\left(\frac{1}{a_n}\int_{(0,\delta\land a_n]}Lf\left(X^x_1(t-)\right)\dd t\right).
\end{align*}
Since,
\[
\min_{y\in[x,\,x+ca_n]}Lf(y)\mathbb E\left(\frac{\delta\land a_n}{a_n}\right)
\le \mathbb E\left(\frac{1}{a_n}\int_{(0,\delta\land a_n]}Lf\left(X^x_1(t-)\right)\dd t\right)
\le\max_{y\in[x,\,x+ca_n]}Lf(y)\mathbb E\left(\frac{\delta\land a_n}{a_n}\right),
\]
we thus infer by means of Lemma~\ref{l.p.t.1.2.4} and the DCT that
\begin{equation}\label{e.jumps_2}
\lim_{n\to\infty}\frac{\mathbb E\left(\sum_{i\in\mathcal I}\mathds1_{(0,\delta\land a_n]}(t_i)\Delta Z^{x,f}_{t_i}\right)}{a_n}
=Lf(x).
\end{equation}

Let us now deal with the remaining contribution to the difference $Z^{x,f}_{a_n}-Z^{x,f}_{0}$. For this purpose,  consider the process $(\hat Z^{x,f}_t)_{t\in\rpn}$, given by
\[
\hat Z^{x,f}_t:=Z^{x,f}_t-\sum_{i\in\mathcal I:t_i\le t}\Delta Z^{x,f}_{t_i}.
\]
Since, according to Lemma~\ref{l.MVT} and (\ref{finite_derivative}),
\[
\mathbb E\left(\sup_{n\in\N}\left|\frac{f(x+c(\delta\land a_n))-f(x)}{c(\delta\land a_n)}\cdot\frac{\delta\land a_n}{a_n}\right|\right)
\le\mathbb E\left(\sup_{y\in(x,x+c(\delta\land a_n))}|f'_+(y)|\right)
\le\sup_{y\in(x,x+c)}|f'_+(y)|<\infty,
\]
we deduce by applying the DCT that
\begin{equation}\label{e.derivative_f}
\begin{aligned}
 \lim_{n\to\infty}\frac{\mathbb E\left(\hat Z^{x,f}_{\delta\land a_n}-\hat Z^{x,f}_0\right)}{a_n}
&= \lim_{n\to\infty}\mathbb E\left(\frac{f(x+c(\delta\land a_n))-f(x)}{a_n}\right)
\\[0.5ex]
&= c\,\lim_{n\to\infty}\mathbb E\left(\frac{f(x+c(\delta\land a_n))-f(x)}{c(\delta\land a_n)}\cdot\frac{\delta\land a_n}{a_n}\right)
\\[0.5ex]
&= c\,\mathbb E\left(\lim_{n\to\infty}\frac{f(x+c(\delta\land a_n))-f(x)}{c(\delta\land a_n)}\cdot\lim_{n\to\infty}\frac{\delta\land a_n}{a_n}\right)
\\[0.5ex]
&=cf'(x).
\end{aligned}
\end{equation}

Combining (\ref{e.jumps_2}) with (\ref{e.derivative_f}) yields that
\begin{align*}
0 &=\lim_{n\to\infty}\frac{\mathbb E\left(Z^{x,f}_{\delta\land a_n}-Z^{x,f}_{0}\right)}{a_n}\notag
\\[0.5ex]
&= \lim_{n\to\infty}\frac{\mathbb E\left(\sum_{i\in\mathcal I}\mathds1_{\{t_i\in(0,\delta\land a_n]\}}\Delta Z^{x,f}_{t_i}\right)}{a_n}+\lim_{n\to\infty}\frac{\mathbb E\left(\hat Z^{x,f}_{\delta\land a_n}-\hat Z^{x,f}_0\right)}{a_n}\notag
\\[0.5ex]
&= Lf(x)+cf'(x)
\end{align*} 
holds for all $x\in\mathcal C_f$, where the first equality results from the martingale property of $Z^{x,f}$ in conjunction with the optional sampling theorem. Consequently, $f$ solves (\ref{e.FKPP.1}), which completes the proof.
\hfill$\square$

\section{Analytic properties of one-sided FKPP travelling waves}\label{s.ap}
In this section we provide the proof of Theorem~\ref{t.ap}. For this purpose we shall resort to the following version of the fundamental theorem of calculus for Dini derivatives, taken from \cite{HT06}.
\begin{proposition}[Theorem~11 in \cite{HT06}]\label{fTCDD}
Let  $a,b\in\R$ with $a<b$. If f is a continuous function that has a finite Dini derivative $f'_+(y)$ for every $y\in[a,b]$, then
\begin{equation}\label{e.fTCDD}
f(b)-f(a) = \int_{[a,b]} f'_+(y)\d y,
\end{equation}
provided that $f'_+$ is Lebesgue integrable over $[a, b]$. 
\end{proposition}

This version of the fundamental theorem of calculus for Dini derivatives will be used in  the proof of Proposition~\ref{l.uniqueness} that we are now going to present. Furthermore, we shall resort to Proposition~\ref{fTCDD} also in the proof of Theorem~\ref{t.ap}, where we show that travelling waves are continuously differentiable.
Let us point out that $f$ having  finite Dini derivatives is essential in Proposition~\ref{fTCDD}. Indeed, for singular functions $f$, such as the Cantor function, the equality in (\ref{e.fTCDD}) does not hold true, since in that case $f'_+=0$ Lebesgue-a.e. but $f$ is not a constant function.

Let us proceed with the following proposition that shows uniqueness of one-sided FKPP travelling waves in $\mathcal T$ with wave speed $c>c_{\bar p}$. Our method of proof for this result makes use of Proposition~\ref{fTCDD}.

\begin{proposition}\label{l.uniqueness}
Any one-sided FKPP travelling wave $f\in\mathcal T$ with wave speed $c>c_{\bar p}$ satisfies 
\[
f=\varphi.
\]
\end{proposition} 

\begin{pf}
Let $f\in\mathcal T$ be a  function that solves (\ref{e.FKPP.1}) and fix some $x>0$. 
Notice first that the map
\[
t\mapsto  \hat Z^{x,f}_t:=Z^{x,f}_t-\sum_{i\in\mathcal I:t_i\le t}\Delta Z^{x,f}_{t_i}
\]
is continuous. In addition, recall from (\ref{d.delta}) the definition
\[
\delta:=\inf\left\{t>0:x+ct+\ln\left(|\pi(t)|^\downarrow_2\right)>0\right\}\land\tau^-_{1,0}\land1
\]
and observe by means of (\ref{finite_derivative}) that the  Dini derivative $\hat Z^{x,f}_+$, given by
\begin{equation}\label{e.l.t.1.1.0.1}
\hat Z^{x,f}_+(s):=\limsup_{h\downarrow0}\frac{\hat Z^{x,f}_{s+h}-\hat Z^{x,f}_s}{h} =c f'_+(X^x_1(s))
\end{equation}
for all $s\in[0,\delta]$, is a finite Lebesgue measurable function, since $f'_+$ and $s\mapsto X^x_1(s)$ are Lebesgue measurable.
Moreover, 
in view of  (\ref{finite_derivative}) we also infer that
\[
\int_{[0,\delta]}\left|\hat Z^{x,f}_+(s)\right|\dd s
\le c\int_{[0,\delta)}\left|f'_+(X^x_1(s))\right|\dd s
\le c \int_{[0,1)}\sup_{y\in[x,\,x+c]}\left|f'_+(y)\right|\dd s
= c\sup_{y\in[x,\,x+c]}\left|f'_+(x)\right|
<\infty
\]
holds $\mathbb P$-almost surely. Hence, $\hat Z^{x,f}_+(s)$ is Lebesgue integrable over $[0,\delta]$. According to Proposition~\ref{fTCDD} we thus obtain   that
\begin{equation}\label{e.l.t.1.1.0.2}
Z^{x,f}_{\delta\land a_n}-Z^{x,f}_0=\hat Z^{x,f}_{\delta\land a_n}-\hat Z^{x,f}_0+\sum_{i\in\mathcal I:t_i\le \delta\land a_n}\Delta Z^{x,f}_{t_i}=\int_{[0,\delta\land a_n]} \hat Z^{x,f}_+(s)\d s+\sum_{i\in\mathcal I:t_i\le \delta\land a_n}\Delta Z^{x,f}_{t_i}
\end{equation}
for all $n\in\N$, where $(a_n)_{n\in\N}$ is a sequence in $(0,1)$ with $a_n\downarrow 0$ as $n\to\infty$.

With $\eta$ being the Poisson random measure on $\rpn\otimes\mathcal P$ that determines $X_1$, we deduce from (\ref{e.l.t.1.1.0.2}), in conjunction with  the compensation formula for Poisson point processes and Fubini's theorem that
\begin{align}\label{e.l.t.1.1.1.2b}
&\mathbb E\left(Z^{x,f}_{\delta\land a_n}\right)-\mathbb E\left(Z^{x,f}_0\right)\notag
\\[0.5ex]
&= \mathbb E\left(\int\limits_{[0,\delta\land a_n]} \hat Z^{x,f}_+(s)\d s\right)\notag
\\
&\qquad +\mathbb E\left(\int\limits_{[0,1]\times\mathcal P}\mathds1_{[0,\delta\land a_n]}(s)\left(\prod_{k\in\N}f(X^x_1(s-)+\ln(|\pi_k|))-f(X^x_1(s-))\right)\eta(\dd s,\dd\pi)\right)
\\[0.5ex]
&= \mathbb E\left(\int\limits_{[0,\delta\land a_n]} \hat Z^{x,f}_+(s)\d s\right)+\mathbb E\left(\int\limits_{[0,\delta\land a_n]}\int\limits_{\mathcal P}\left(\prod_{k\in\N}f(X^x_1(s-)+\ln(|\pi_k|))-f(X^x_1(s-))\right)\mu(\dd\pi)\d s\right)\notag
\\[0.5ex]
&= \mathbb E\left(\int\limits_{[0,\delta\land a_n]} \hat Z^{x,f}_+(s)\d s\right)+\mathbb E\left(\int\limits_{[0,\delta\land a_n]}Lf(X^x_1(s-)\d s\right)\notag
\end{align}
for all $n\in\N$.
Observe that $f$ being monotone and $X^x_1$ having only at most countably many jumps $\mathbb P$-a.s. implies that $X^x_1(s)\in C_f$ for Lebesgue-a.a. $s\in(0,1)$ $\mathbb P$-almost surely. By means of
(\ref{e.l.t.1.1.0.1}) and (\ref{e.l.t.1.1.1.2b}) as well as the fact that any $u\in(0,1)$ is $\mathbb P$-a.s. not a jump time of $\Pi$ this results in
\[
\mathbb E\left(Z^{x,f}_{\delta\land a_n}\right)-\mathbb E\left(Z^{x,f}_0\right)
= \mathbb E\left(\int_{[0,\delta\land a_n]}\left(c f'_++Lf\right)(X^x_i(s))\dd s\right)= 0,
\]
i.e.
\[
\mathbb E\left(Z^{x,f}_{\delta\land a_n}\right)=\mathbb E\left(Z^{x,f}_0\right)=f(x).
\]
By means of the strong fragmentation property of $\Pi$ we thus conclude that
\[
\mathbb E\left(\left.Z^{x,f}_{\tau+\delta\land a_n}\right|\mathscr F_\tau\right)=\prod_{n\in\mathcal N^x_\tau}\left.\mathbb E\left(Z^{y,f}_{\delta\land a_n}\right)\right|_{y=X^x_n(\tau)}=\prod_{n\in\mathcal N^x_\tau}f(X^x_n(\tau))=Z^{x,f}_\tau
\]
holds $\mathbb P$-a.s. for every 
finite stopping time $\tau$. Therefore,
\begin{equation}\label{uniqueness_martingale}
\begin{aligned}
\mathbb E\left(\left.Z^{x,f}_{t+k(\delta\land a_n)}\right|\mathscr F_t\right)
&= \mathbb E\left(\left.Z^{x,f}_{t+\delta}+\sum_{n=2}^k\mathbb E\left(\left.Z^{x,f}_{t+n\delta}-Z^{x,f}_{t+(n-1)\delta}\right|\mathscr F_{t+(n-1)\delta}\right)\right|\mathscr F_t\right)
\\
&= \mathbb E\left(\left.Z^{x,f}_{t+\delta}\right|\mathscr F_t\right)
\\
&= Z^{x,f}_t
\end{aligned}
\end{equation}
$\mathbb P$-a.s. for all $t\in\rpn$ and every measurable $k:\Omega\to\N$. For any $s,t\in\rpn$ set
\[
k_s:=\left\lfloor\frac{s}{\delta\land a_n}\right\rfloor\qquad\text{as well as}\qquad r_s:=\frac{s}{\delta\land a_n}-k_s\in(0,1)
\]
and observe that $Z^{x,f}$ is $\mathbb P$-a.s. left-continuous at $t+s$, since $t+s$ is $\mathbb P$-a.s.  not a jump time of $\Pi$. In conjunction with the DCT for conditional expectations and (\ref{uniqueness_martingale}) this implies that
\[
\mathbb E\left(\left.Z^{x,f}_{t+s}\right|\mathscr F_t\right)
=\lim_{n\to\infty}\mathbb E\left(\left.Z^{x,f}_{t+s-r_s(\delta\land a_n)}\right|\mathscr F_t\right)
=\lim_{n\to\infty}\mathbb E\left(\left.Z^{x,f}_{t+k_s(\delta\land a_n)}\right|\mathscr F_t\right)=Z^{x,f}_t
\]
for all $s,t\in\rpn$.
Hence, $Z^{x,f}$ is a martingale and thus we deduce from Theorem~10 in \cite{KK12} that $f=\varphi$.
\end{pf}

Proposition~\ref{l.uniqueness} shows that in order to derive analytic properties of one-sided FKPP travelling waves in $\mathcal T$ we only need to consider the function $\varphi$. Bearing this in mind we proceed to prove Theorem~\ref{t.ap}. In order to obtain strict monotonicity of $\varphi$ we shall use the following result.

\begin{lemma}\label{l.continuity.1.a0}
Let $c>c_{\bar p}$.  For any $0\le x<y<\infty$ there exists some $\alpha_{x,y}>0$\label{p.a_x_y} such that
\[
\varphi(x)-\varphi(x+h)\ge\alpha_{x,y}\left(\varphi(y)-\varphi(y+h)\right)
\]
for all $h>0$.
\end{lemma}

\begin{pf}
In the first part of this proof we show that for every deterministic time $t>0$ the probability that $X_1$ reaches level $x>0$ before time $t$ is positive. In the second part we use this fact in order to obtain a lower bound of the probability that  for some $n\in\N$ the process $X^x_n$ hits a given level $y>x$ before some deterministic time $s>0$. Subsequently, we combine this lower bound with the estimate  (\ref{e.N_x}) of the number of blocks that are alive at a given time and with the positivity of the probability of extinction.

\underline{Part I}
For every $x\in\rpn$ set
\[
\quad \tau^+_{1,x}:=\inf\{t\in\rpn:X_1(t)>x\}.
\]
According to Corollary~3.14 in \cite{103} we have that $(\tau^+_{1,x})_{x\in\rpn}$ is a subordinator with either killing at an independent exponential ``time'' $\mathfrak e$ or  with no killing in which case we set $\mathfrak e:=\infty$. Moreover, by means of Proposition~1.7 in \cite{Ber99} we thus infer that
\begin{equation}\label{e.estimate.partI.b}
\mathbb P\left(\tau^+_{1,x}< t\right)=\mathbb P\left(\{\tilde\tau^+_{1,x}< t\}\cap\{x<\mathfrak e\}\right)=\mathbb P\left(\tilde\tau^+_{1,x}< t\right)\mathbb P\left(x< \mathfrak e\right)>0
\end{equation}
holds for all $t>0$ and $x\in\rpn$, where $(\tilde\tau^+_{1,x})_{x\in\rpn}$ is some non-killed subordinator, independent of $\mathfrak e$, satisfying 
\[
\tilde\tau^+_{1,x}\mathds1_{\{x< \mathfrak e\}}=\tau^+_{1,x}\mathds1_{\{x< \mathfrak e\}}.
\] 
For the time being, fix some $x\in\rpn$. Let us now show that 
\begin{equation}\label{e.estimate.partI.b2}
\forall\,t>0:\,\mathbb P(\tau^+_{1,x}<\tau^-_{1,0}\land t)>0.
\end{equation} 
To this end, assume  we have
\begin{equation}\label{e.estimate.partI.b2b}
\exists\,t_0>0:\,\mathbb P(\tau^+_{1,x}<\tau^-_{1,0}\land t_0)=0.
\end{equation}
Our goal is to show that this results in a contradiction. For this purpose, set $\tau^2_0:=\tilde\tau^2_0:=0$ and for every $n\in\N$ define
\[
\tilde\tau^1_n :=  \inf\{t\ge\tilde\tau^2_{n-1}:X_1(t)<0\}
\qquad\text{as well as}\qquad
\tilde\tau^2_n := \inf\{t\ge\tilde\tau^1_n:X_1(t)=0\}.
\]
In addition, set
\[
n^*:=\sup\left\{n\in\N:\tilde\tau^2_n<\infty\right\}
\]
as well as
\[
\tau^1_n :=  \inf\{t\ge\tau^2_{n-1}:X_1(t)<0\}
\qquad\text{and}\qquad
\tau^2_n := \inf\{t\ge\tau^1_n:X_1(t)=0\}\land \tilde\tau^2_{n^*},
\]
where $\tilde\tau^2_\infty:=\infty$.
Since  for $X_1$ the point $0$ is irregular for $(-\infty,0)$, there exists some $\varepsilon>0$ such that $\mathbb P(\tau^-_{1,0}\ge\varepsilon)>0$ and consequently we obtain by means of the strong Markov property of $\Pi$ that
\begin{equation}\label{e.c.l.l.p.2.2.0.1b}
\sum_{n\in\N}\mathbb P\left(\tau^1_n-\tau^2_{n-1}\ge\varepsilon\left|\mathscr F_{\tau^2_{n-1}}\right.\right)=\sum_{n\in\N}\mathbb P\left(\tau^-_{1,0}\ge\varepsilon\right)=\infty
\end{equation}
$\mathbb P$-almost surely. Since $\{\tau^1_n-\tau^2_{n-1}\ge\varepsilon\}$ is $\mathscr F_{\tau^2_{n}}$-measurable, we can apply an extended Borel-Cantelli lemma (see e.g. \cite[(3.2) Corollary in Chapter~4]{Dur91} or \cite[Corollary~5.29]{113}) to deduce that
\[
\left\{\{\tau^1_n-\tau^2_{n-1}\ge\varepsilon\}\text{ holds for infinitely many $n\in\N$}\right\}=\left\{\sum_{n\in\N}\mathbb P\left(\tau^1_n-\tau^2_{n-1}\ge\varepsilon\left|\mathscr F_{\tau^2_{n-1}}\right.\right)=\infty\right\}.
\]
Thus (\ref{e.c.l.l.p.2.2.0.1b}) implies that $\tau^2_n\to\infty$ $\mathbb P$-a.s. on the event $\{n^*=\infty\}$ as $n\to\infty$. 
With $t_0$ given by (\ref{e.estimate.partI.b2b}) another application of the strong Markov property therefore yields that
\begin{equation}\label{e.estimate.partI.c}
\mathbb P\left(\tau^+_{1,x}< t_0\right) \le \mathbb E\left(\sum_{n\in\N}\mathbb P\left(\tau^3_{n,x}<\tau^1_n\land t_0\left|\mathscr F_{\tau^2_{n-1}}\right.\right)\right)
= \sum_{n\in\N}\mathbb P\left(\tau^+_{1,x}<(\tau^-_{1,0}\land t_0)\right)
 =0,
\end{equation}
where 
\[
\tau^3_{n,x} :=  \inf\{t\ge\tau^2_{n-1}:X_1(t)>x\}
\]
for all $n\in\N$. Since (\ref{e.estimate.partI.c}) contradicts (\ref{e.estimate.partI.b}), we conclude that (\ref{e.estimate.partI.b2}) does indeed hold true.

\underline{Part II}
Let $0\le x<y<\infty$ and  for any $t\in\rpn$ set $R^x_1(t):=\sup_{n\in\N}X^x_n(t)$. In addition, we define
\[
\tau^+_y(x):=\inf\left\{t\in\rpn:R^x_1(t)\ge y\right\}. 
\]
Note that $R^x_1(\tau^+_y(x))=y$ if $\tau^+_y(x)<\infty$, since $R^x_1$ does not jump upwards and thus creeps over the level $y$. Furthermore, let $s>0$ and set $\gamma:=e^{x+cs}-1$ as well as
\[
\alpha_{x,y}:=\mathbb P\left(\tau^+_y(x)<\zeta^x\land s\right)\mathbb P(\zeta^{y}<\infty)^\gamma.
\]
Observe that (\ref{e.estimate.partI.b2}) and Proposition~\ref{positivesurviaval} imply that $\alpha_{x,y}>0$,
since 
\[
\mathbb P\left(\tau^+_y(x)<\zeta^x\land s\right)\ge\mathbb P\left(\tau^+_{1,y-x}<\tau^-_{1,0}\land s\right).
\] 

By means of the strong fragmentation property of $\Pi$ we deduce that
\begin{align*}
\varphi(x)-\varphi(x+h) &= \mathbb P(\zeta^x<\infty)-\mathbb P(\zeta^{x+h}<\infty)\notag
\\[0.5ex]
&\overset{(*)}\ge \mathbb P\left(\tau^+_y(x)<\zeta^x\land s\right)\mathbb P(\zeta^y<\infty)^\gamma\left(\mathbb P(\zeta^y<\infty)-\mathbb P(\zeta^{y+h}<\infty)\right)
\\[0.5ex]
&= \alpha_{x,y}\left(\varphi(y)-\varphi(y+h)\right)\notag
\end{align*}
holds true for any $h>0$, where the exponent $\gamma$ in $(*)$ results from the estimate
\[
N^x_{\tau^+_y(x)}\le e^{x+c\tau^+_y(x)}< e^{x+cs}=\gamma+1
\]
$\mathbb P$-a.s. on $\{\tau^+_y(x)<s\}$. Notice that in $(*)$ we have used that the value of $X^x_n$, $n\in\N$, at time $\tau^+_y(x)$ is less than or equal to $y$ as well as the monotonicity of the  probability of extinction.
\end{pf}

Observe that $\varphi$ is clearly a monotone function. However, even though monotonicity is trivial, it is not obvious whether  $\varphi$ is strictly monotone. The following lemma answers the question regarding strict monotonicity of $\varphi$ affirmatively.

\begin{lemma}\label{positivesurviaval.monotonicity}
Let $c>c_{\bar p}$. Then $\varphi$ is strictly monotonically decreasing on $\rpn$.
\end{lemma}

\begin{pf}
Let $x\in\rpn$ and set
\[
\gamma_x:=\ln\left(|\pi(\zeta^x)|^\downarrow_1\cdot\left|\Pi^x_{\kappa(\zeta^x)}(\zeta^x-)\right|\right).
\]
According to Proposition~\ref{positivesurviaval} we have $\mathbb P(\zeta^x<\infty)>0$ and hence 
\begin{align*}
& \mathbb P\left(\{\zeta^x<\infty\}\cap\bigcup_{n\in\N}\{x+c\zeta^x+\gamma_x\in(-n,0)\}\right) 
\\[0.5ex]
&= \mathbb P(\{\zeta^x<\infty\}\cap\{x+c\zeta^x+\gamma_x\in (-\infty,0)\})
\\[0.5ex]
&= \mathbb P(\zeta^x<\infty)
\\[0.5ex]
&> 0.
\end{align*}
Therefore,  there exists some $z>0$ such that
\[
\mathbb P\left(\{\zeta^x<\infty\}\cap\{x+c\zeta^x+\gamma_x\in(-z,0)\}\right)>0
\]
and thus the strong fragmentation property, in conjunction with Proposition~\ref{positivesurviaval}, yields that
\[
\mathbb P(\{\zeta^x<\infty\}\cap\{\zeta^{x+z}=\infty\})\ge \mathbb P(\{\zeta^x<\infty\}\cap\{x+c\zeta^x+\gamma_x\in(-z,0)\})\mathbb P(\zeta^0=\infty)
> 0.
\]
Consequently,  there exists some $z>0$ such that
\begin{align}\label{e.l.positivesurviaval.2.1c5}
\mathbb P(\zeta^x<\infty) &= \mathbb P(\{\zeta^x<\infty\}\cap\{\zeta^{x+z}=\infty\})+\mathbb P(\{\zeta^x<\infty\}\cap\{\zeta^{x+z}<\infty\})\notag
\\[0.5ex]
&> \mathbb P(\zeta^{x+z}<\infty),
\end{align}
where the final estimate follows from $\{\zeta^{x+z}<\infty\}\subseteq\{\zeta^x<\infty\}$. 

Observe that (\ref{e.l.positivesurviaval.2.1c5}) implies that for every $h>0$  there exists some $y\ge x$ such that 
\begin{equation}\label{e.l.positivesurviaval.2.1e}
\varphi(y)>\varphi(y+h).
\end{equation} 

According to Lemma~\ref{l.continuity.1.a0}, for all $h>0$ and $y\ge x$ satisfying (\ref{e.l.positivesurviaval.2.1e}) there exists some $\alpha_{x,y}>0$ such that
\[
\varphi(x)-\varphi(x+h)
\ge \alpha_{x,y}\left(\varphi(y)-\varphi(y+h)\right)
>0.
\]
Since $x\in\rpn$ was chosen arbitrarily, this proves the assertion that $\varphi$ is strictly monotonically decreasing on $\rpn$.
\end{pf}

In the proof of Theorem~\ref{t.ap} we shall make use of the theory of scale functions for spectrally negative L\'evy processes. For this purpose, let $W$ be the scale function of the spectrally negative L\'evy process $X_1$. That is to say, $W$ is the unique continuous and strictly monotonically increasing function $W:\rpn\to\rpn$, whose Laplace transform satisfies
\[
\int_{(0,\infty)}e^{-\beta x}W(x)\d x=\frac{1}{\psi(\beta)}
\]
for all $\beta>\Psi(0)$, 
where $\psi$ denotes the Laplace exponent of $X_1$ and $\Psi(0):=\sup\{\lambda>0:\psi(\lambda)=0\}$.

Let us now tackle the proof of Theorem~\ref{t.ap}.

{\bf Proof of Theorem~\ref{t.ap}}
This proof is divided into two parts. In the first part  we show that $\varphi\in \mathcal T$ and that $\varphi$ is right-continuous at $0$. Subsequently, in the second part we use the continuity of $\varphi|_{\rp}$  in order to prove that $\varphi|_{\rp}$ is continuously differentiable, if $\varphi$ solves (\ref{e.FKPP.1}). According to Proposition~\ref{l.uniqueness} and Lemma~\ref{positivesurviaval.monotonicity} the proof is then complete.

\underline{Part I}
Recall the definition of $\mathcal T$ in Definition~\ref{d.T} and note that  Theorem~10 in \cite{KK12} yields that $\varphi$ satisfies (\ref{e.FKPP.2}). Hence, since $\varphi$ is nonincreasing, in order to prove $\varphi\in\mathcal T$ it remains to show that $\varphi|_{\rp}$ is continuous and that (\ref{finite_derivative}) holds.  To this end, let  $\mathfrak n$ denote the excursion measure of excursions $({\bf e}_s)_{s\in\rpn}$ indexed by the local time at the running maximum of the L\'evy process $X_1$. Furthermore, for any such excursion ${\bf e}$ let $\bar{\bf e}$ denote the height  of this excursion.
According to Lemma~8.2 in \cite{103} the scale function $W$ has a right-derivative on $\rp$ given by
\[
W'_+(x)=W(x)\mathfrak n(\bar{\bf e}> x) 
\]
for any $x\in\rp$.  Note that $\mathfrak n$ being $\sigma$-finite (cf. Theorem~6.15 in \cite{103})   implies that 
\[
\sup_{y\ge x}\mathfrak n(\bar{\bf e}> y)=\mathfrak n(\bar{\bf e}> x)<\infty.
\]
for every $x\in\rpn$. Moreover, we have
\[
\mathbb P(\xi^x=\infty)\ge\mathbb P_x\left(\tau^+_{1,y}<\tau^-_{1,0}\right)\mathbb P(\xi^y=\infty)
\]
for all $x,y\in\rpn$ with $x<y$, where under $\mathbb P_x$ the process $X_1$ is shifted to start in $x$. Therefore, we obtain
\begin{equation}\label{e.t.mainresult.1.2.p2}
\varphi(x)-\varphi(x+h)\le(1-\varphi(x))\left(\frac{1}{\mathbb P_x\left(\tau^+_{1,x+h}<\tau^-_{1,0}\right)}-1\right)\le\frac{1}{\mathbb P_x\left(\tau^+_{1,x+h}<\tau^-_{1,0}\right)}-1
\end{equation}
for all $h\in\rp$ and $x\in\rpn$. By means of (8.8) in Theorem~8.1 of \cite{103} we have
\begin{equation}\label{e.t.mainresult.1.2.p2_c}
\frac{1}{\mathbb P_x\left(\tau^+_{1,x+h}<\tau^-_{1,0}\right)}-1=\frac{W(x+h)-W(x)}{W(x)}
\end{equation}
for any $x\in\rpn$. Consequently,  
\[
\sup_{y\ge x}\left|\varphi'_+(y)\right|\le\sup_{y\ge x}\limsup_{h\downarrow0}\left(\frac{1}{W(y)}\frac{W(y+h)-W(y)}{h}\right)=\sup_{y\ge x}\mathfrak n(\bar{\bf e}> y)<\infty
\] 
holds for every $x\in\rp$. Moreover, in view of (\ref{e.t.mainresult.1.2.p2}), (\ref{e.t.mainresult.1.2.p2_c}) and the continuity of $W|_{\rpn}$ we deduce  that $\varphi|_{\rp}$ is continuous. Therefore, we conclude that $\varphi\in\mathcal T$. Since $X_1$ has bounded variation, we infer by means of  Lemma~8.6 in \cite{103} that $W(0)>0$ and thus the above line of argument also  yields  that $\varphi$ is right-continuous at $0$.

\underline{Part II}
Assume that $\varphi$ satisfies (\ref{e.FKPP.1}) on $\mathcal C_\varphi$. In view of Part~I it follows from Lemma~\ref{l.p.t.1.2.4} that $L\varphi$ is continuous. Moreover, we deduce from  (\ref{e.FKPP.1}) and the monotonicity of $\varphi$ that $\varphi'_+=-c^{-1}L\varphi$ Lebesgue-almost everywhere on $\rp$. Since  the upper Dini derivative $\varphi'_+$ is bounded on any interval $[a,b]\subseteq\rp$,  it thus follows from Proposition~\ref{fTCDD} and Lebesgue's integrability criterion for Riemann integrals that
\[
\varphi(b)-\varphi(a)=\int_a^b\varphi'_+(x)\d x=-\frac{1}{c}\int_a^bL\varphi(x)\d x=F(b)-F(a)
\]
for Lebesgue-almost all $a,b\in\rp$, where  $F\in C^1(\rp,\R)$ is an antiderivative of $-c^{-1}L\varphi$ on $\rp$. Hence, we have $\varphi=F+\text{const.}$ Lebesgue-almost everywhere on $\rp$. Since $\varphi|_{\rpn}$ and $F$ are continuous, this implies that $\varphi|_{\rp}=F+\text{const.}$ 
and consequently  $\varphi|_{\rp}\in C^1(\rp,[0,1])$.
In the light of  Proposition~\ref{l.uniqueness} and Lemma~\ref{positivesurviaval.monotonicity} this proves the assertion.
\hfill$\square$

\section{Existence and uniqueness of one-sided travelling waves}\label{s.mainresult.1}

This section is devoted to the proof of Theorem~\ref{t.mainresult.1}. Our method of proof makes use of the results that we developed in the previous two sections.

{\bf Proof of Theorem~\ref{t.mainresult.1}}
The first part of the proof shows the nonexistence of one-sided FKPP travelling waves in $\mathcal T$ for wave speeds $c\le c_{\bar p}$ and the second part proves the existence  of such travelling waves for wave speeds above the critical value $c_{\bar p}$. The uniqueness was shown in Proposition~\ref{l.uniqueness}. 

\underline{Part I}
Fix some $c\le c_{\bar p}$ as well as $x\in\rpn$ and let $f\in\mathcal T$. Further, assume that $f$ satisfies (\ref{e.FKPP.1}). Then the proof of Proposition~\ref{l.uniqueness} shows that $(Z^{x,f}_t)_{t\in\rpn}$ is a uniformly integrable martingale and hence the $\mathbb P$--a.s. martingale limit $Z^{x,f}_\infty:=\lim_{t\to\infty}Z^{x,f}_t$ satisfies
\begin{equation}\label{e.l.t.mainresult.1.1.1}
\mathbb E\left(Z^{x,f}_\infty\right)=\mathbb E\left(Z^{x,f}_0\right)=f(x).
\end{equation}
Since $c\le c_{\bar p}$, we have according to Proposition~\ref{positivesurviaval} that $\mathbb P(\zeta^x<\infty)=1$, that is to say $\mathcal N^x_t\to\emptyset$ $\mathbb P$-a.s. as $t\to\infty$. Because the empty product equals 1, we thus infer that
\[
Z^{x,f}_\infty=\lim_{t\to\infty}\prod_{n\in\mathcal N^x_t}f(X^x_n(t))=1
\]
$\mathbb P$-almost surely. In view of (\ref{e.l.t.mainresult.1.1.1}) this implies that $f\equiv1$, which is a contradiction to $f\in\mathcal T$, since every $f\in\mathcal T$ satisfies (\ref{e.FKPP.2}). Consequently, there does not exist a 
function $f\in\mathcal T$ that satisfies (\ref{e.FKPP.1}). 

\underline{Part II}
Now let  $c>c_{\bar p}$ and $x\in\rpn$. In the light of Proposition~\ref{l.uniqueness} it only remains to show that $\varphi$ is indeed a one-sided FKPP travelling wave with wave speed $c$. Since $\varphi\in\mathcal T$ satisfies the boundary condition (\ref{e.FKPP.2}), we only have to deal with (\ref{e.FKPP.1}). In order to prove that $\varphi$ solves  (\ref{e.FKPP.1}) we aim at applying Theorem~\ref{p.t.1.2}. To this end, observe that the fragmentation property of $\Pi$ yields that
\[
\varphi(x)=\mathbb E(\mathbb P(\zeta^x<\infty|\mathscr F_t))=\mathbb E\left(\prod_{n\in\mathcal N^x_t}\varphi(X^x_n(t))\right)=\mathbb E\left(Z^{x,\varphi}_t\right)
\]
for every $t\in\rpn$. By means of another application of the fragmentation property  we therefore deduce that
\[
\mathbb E\left(\left.Z^{x,\varphi}_{t+s}\right|\mathscr F_t\right)=\prod_{n\in\mathcal N^x_t}\left.\mathbb E\left(Z^{y,\varphi}_s\right)\right|_{y=X^x_n(t)}=\prod_{n\in\mathcal N^x_t}\varphi(X^x_n(t))=Z^{x,\varphi}_t
\]
holds $\mathbb P$-a.s. for all $s,t\in\rpn$. Hence, $Z^{x,\varphi}$ is a $\mathbb P$-martingale. In the proof of Theorem~\ref{t.ap} we have shown that  $\varphi\in \mathcal T$ and consequently we infer from Theorem~\ref{p.t.1.2} that $\varphi$ solves the integro-differential equation (\ref{e.FKPP.1}). 
\hfill$\square$


\end{document}